%% file: submission.2.tex
\documentclass[twocolumn]{autart}    

\input{headerSetupAutomatica.tex}

\begin{document}

\begin{frontmatter}

\title{On Physical Realizability of 
	Nonlinear Quantum Stochastic Differential Equations} 

\thanks[footnoteinfo]{This work was supported by the
	Australian Research Council under grant numbers FL110100020 and DP180101805. It was also supported by the Airforce Office Scientific Research (AFOSR) under agreement number FA2386-16-1-4065.}

%
\author[ADFA,ANU]{Muhammad F. Emzir}\ead{puat133@gmail.com},    
\author[ADFA]{Matthew J. Woolley}\ead{m.woolley@adfa.edu.au},               
\author[ANU]{Ian R. Petersen}\ead{i.r.petersen@gmail.com}  
\address[ADFA]{School of Engineering and Information Technology, University of New South Wales, ADFA, Canberra, ACT 2600, Australia}
\address[ANU]{Research School of Engineering, Australian National University, Canberra ACT 0200, Australia}


\begin{abstract}                          
In this article we study \emph{physical realizability} for a class of nonlinear quantum stochastic differential equations (QSDEs). Physical realizability is a property in which a QSDE corresponds to the dynamics of an open quantum system. We derive a sufficient and necessary condition for a nonlinear QSDE to be physically realizable. 
\end{abstract}

\end{frontmatter}

\section{Introduction}
The theory of coherent quantum feedback control has attracted a significant amount of interest \cite{Wiseman1994,Yanagisawa2003,Yanagisawa2003a,James2008,Nurdin2009}. A coherent feedback controller is a quantum system that is constructed to coherently manipulate the output field of the controlled system and feed it back to the controlled system. Another approach to implement feedback control for a quantum system is based on measurements of the quantum system, where the control input is computed based on the measurement record. This approach, known as measurement-based feedback control (MBFC), has been well studied over the last two decades \cite{Wiseman1993,Wiseman1994a,Belavkin1999,wiseman2010quantum}. 
\\
There are situations in which coherent feedback control potentially offers advantages over MBFC; e.g, see \cite{Nurdin2009,Hamerly2012,Yamamoto2014}. 
Coherent feedback has been demonstrated in  optics \cite{Mabuchi2008,Iida2012,Zhou2015}, superconducting circuits \cite{Kerckhoff2012}, and electromechanical systems \cite{Kerckhoff2013}.
\\
In designing a coherent quantum feedback controller, the question arises as to whether a given feedback controller can be realized as a physical quantum system.  A Markovian open quantum system's dynamics can be characterized by its Hamiltonian operator $\mathbb{H}$, its coupling operator $\mathbb{L}$, and its  scattering matrix $\mathbb{S}$ \cite{Parthasarathy1992}. 
Given these parameters, the dynamics of this open quantum system can be described by quantum stochastic differential equations (QSDEs). In contrast to classical stochastic differential equations, for any given QSDEs, there might not be an open quantum system which has these dynamics. Therefore we consider conditions under which a given QSDE corresponds to the dynamics of an open quantum system. A QSDE which corresponds to the dynamics of an open quantum system is said to be \emph{physically realizable}.
\\
For linear QSDEs, the notion of physical realizability was introduced in \cite{James2008}, which presented some algebraic conditions for physical realizability in terms of the linear dynamics matrices. Furthermore, the papers \cite{Nurdin2009a,Nurdin2010,Petersen2011,Gough2015} have shown how to construct a physical quantum optical system from basic quantum optical components corresponding to a given physically realizable QSDE.
However for nonlinear QSDEs, there are only a few limited results available on physical realizability. The paper \cite{maalouf2012physical} considers physical realizability for a class of nonlinear QSDEs satisfying a strong structural assumption. The paper \cite{Espinosa2016} considers the physical realizability of bilinear QSDEs corresponding to finite level quantum systems.
\\
Physical realizability conditions for nonlinear QSDEs will play an important role in designing nonlinear coherent controllers for quantum systems since any designed controller must be physically realizable \cite{Maalouf2011,Zhang2012}. They will also be useful in the modeling of unknown nonlinearities in linear QSDEs, which helps in the synthesis of robust linear quantum controllers as considered in \cite{Petersen2012,Xiang2015}. Furthermore, the physical realizability of nonlinear QSDEs will lead to constraints which need to be satisfied in the system identification of  nonlinear quantum systems \cite{Guta2016}. 
\\
In this article, we derive a sufficient and necessary condition for a general class nonlinear QSDEs to be physically realizable. We initially consider the single mode case where there are two observables of interest $x_{1}(t),x_{2}(t)$, which satisfy the canonical commutation relation $\commutator{x_{1}(t)}{x_{2}(t)} = i \hbar$. In particular we will consider $x_{1}(t) = q(t)$ the position operator and $x_{2}(t) = p(t)$ the momentum operator.
The dynamics of these observables are assumed to be given by
\begin{dmath}
	d\mathbf{x}(t) = \mathbf{f}(\mathbf{x}(t))dt + \mathbf{g}(\mathbf{x}(t)) dA(t)^\ast + \mathbf{g}(\mathbf{x}(t))^\ast dA(t),\label{eq:HeisenBergNonLinearQSDE}
\end{dmath}
where in this equation, both $\mathbf{f}$ and $\mathbf{g}$ are assumed to be written as a power series of the pair $x_{1}(t),x_{2}(t)$, whilst 
$ dA(t)^\ast,dA(t)$, are the annihilation and creation processes of the quantum field \cite{Parthasarathy1992}. In comparison to the class of nonlinear QSDEs considered in \cite{maalouf2012physical,Espinosa2016}, this class of QSDEs will cover a larger class of nonlinear QSDEs which may arise in many practical applications. Indeed, in many practical applications, both the Hamiltonian and the coupling operator can be modeled by polynomial functions. In addition to the single-mode case, we generalize further the notion of physical realizability and our results to the case of multi-mode quantum systems.
\\
The article is organized as follows. In the Section II, we will review some basic facts on quantum mechanics and open quantum systems. We also introduce new notations about the class of nonlinear QSDEs considered in this article and its properties. The main result for the single mode case will be given in the Section III, where we derive a necessary and sufficient condition for the physical realizability of a QSDE and the preservation of the commutation relations. In essence, we show that both $\mathbf{f}$ and $\mathbf{g}$ have to be conservative vectors of potential operators under suitable choice of axes, which is stated in  Theorem \ref{thm:PhysicalRealizabilityOpen}.
Section IV will cover the generalization of the result obtained in Section III to the multi-mode case. In the Section V, we will discuss some examples of physical realizability of QSDEs. The last section will give the conclusions of this paper.

\subsection{Notation}
We denote the identity operator on a Hilbert space by $1$. Bold letters (e.g. $\mathbf{y}$) will be used to denote a matrix or vector whose elements are Hilbert space operators. Hilbert space adjoints, are indicated by $^{\ast}$, while the adjoint transpose of a vector or matrix of operators will be denoted by $\dagger$; i.e. $\qty(\mathbf{x}^{\ast})^{\top} = \mathbf{x}^{\dagger}$. For single-element operators we will use $*$ and $\dagger$ interchangeably. The commutator of vectors of operators $\mathbf{x}$ and $\mathbf{y}$ is given by $[\mathbf{x}, \mathbf{y}^\top ] = \mathbf{x}\mathbf{y}^\top - \qty(\mathbf{y}\mathbf{x}^\top)^\top$. For a set $A$, $A^{\complement}$ denotes the complement of $A$ with respect to a particular universe. For notation simplicity, we set $\hbar=1$

\section{Preliminaries}
In this section we will present some preliminaries that will be used in the subsequent sections. We will first review the concepts of closed and open quantum systems. Furthermore, in Subsections 2.2-2.4 we will introduce new results about the class of nonlinear functions considered in this article and its properties. In Section \ref{sec:MainResults}, we will build up our main result using these notations. All proofs of the results in these subsections are given in the Appendix.
\subsection{Closed and Open Quantum System Dynamics}
For quantum systems, in contrast to classical systems where the state is determined by a set of scalar variables, the \emph{state} of the system is described by a vector in the system's Hilbert space $\mathscr{h}$. Furthermore, in quantum mechanics, physical quantities like the spin of an atom, position, and momentum, are described as self adjoint operators on a Hilbert space. These operators called observables. The expected values of these quantities are given by an inner product. For example, an observable $A$ at state $\ket{\psi}$ has expected value $\bra{\psi}A\ket{\psi}$ using the Dirac notation; e.g., see \cite{sakurai2013modern}.
\\
The dynamics of a closed quantum system are described by an observable called the \emph{Hamiltonian} $\mathbb{H}$ which acts on the state vectors in the system's Hilbert space, as per
\begin{dmath}
	\dv{\ket{\psi(t)}}{t} = -i \mathbb{H} \ket{\psi(t)},
\end{dmath}
which is  known as the Schr\"{o}dinger equation. The evolution of the state can be described by a unitary operator $U(t)$, where $\ket{\psi(t)} = U(t) \ket{\psi(0)}$. Accordingly, the Schr\"{o}dinger equation can be rewritten as
\begin{dmath}
	{\dv{U(t)}{t} = -i \mathbb{H} U(t), U(0) = 1}.\label{eq:ClosedUnitraryEvolution}
\end{dmath}
From this equation, any system observable $X$ will evolve according to $X(t) = U(t)^\ast X U(t)$, satisfying
\begin{dmath}
	\dv{X(t)}{t} = -i \commutator{X(t)}{\mathbb{H}}, \label{eq:ClosedSystemHeisenberg}
\end{dmath}
which is called the Heisenberg equation of motion for the observable $X$.
\\
An open quantum system is a quantum system which interacts with other quantum mechanical degrees of freedom. Typically, we assume interaction with a large number of degrees of freedom, and refer to them collectively as the environment of our system. An open quantum system $\mathcal{G}$ can be characterized by a triple $\left(\mathbb{S},\mathbb{L},\mathbb{H}\right)$, with  Hamiltonian $\mathbb{H}$, coupling operator $\mathbb{L}$ and  scattering matrix $\mathbb{S}$, which are operators on the \emph{system}'s Hilbert space $\mathscr{h}$. Furthermore, for an open quantum system with $n$ field channels, the matrix $\mathbb{S}$ satisfies, $\mathbb{S}\mathbb{S}^\dagger = \mathbb{S}^\dagger\mathbb{S}=\mathbf{I}$, where $\mathbf{I}$ is the $n \times n$ identity matrix.
\\
The environment is modeled using a symmetric Fock space  $\Gamma\qty(\mathscr{H})$. Specifically,  $\Gamma\qty(\mathscr{H})$ is the Hilbert space for an infinite number of indistinguishable environment particles, where a single environment particle Hilbert space is given by $\mathscr{H}$ \cite{Parthasarathy1992,Reed1972}. 
The total Hilbert space will be given as $\tilde{\mathscr{H}} = \mathscr{h} \otimes \Gamma\qty(\mathscr{H})^{\otimes^n}$.
\\
In a similar way to the unitary operator evolution in the closed quantum system \eqref{eq:ClosedUnitraryEvolution}, we could also derive the unitary operator evolution for an open quantum system. In contrast to the closed quantum system unitary evolution \eqref{eq:ClosedUnitraryEvolution}, the interaction with the environment leads to randomness of the unitary evolution of an open quantum system $\mathcal{G}$ as follows: 
\begin{dmath}
	dU(t) =   \left[\text{tr}\left[\left(\mathbb{S} - \mathbf{I}\right)d\mathbf{\Lambda}^\top(t)\right] + d\mathbf{A}^{\dagger}(t) \mathbb{L}  - \mathbb{L}^{\dagger}\mathbb{S} d\mathbf{A}(t) - \left(\dfrac{1}{2}\mathbb{L}^{\dagger}\mathbb{L}+i\mathbb{H}\right)dt\right]U(t), {\quad U(0) = 1}. \label{eq:UnitaryEvolution}
\end{dmath}
In this equation, $\mathbf{A}(t) = \qty[A_{1}(t) \; \cdots \; A_{n}(t)]$ is a vector of annihilation operators defined on distinct copies of the Fock space $\Gamma\qty(\mathscr{H})$ \cite{nurdin2014quantum}. Each annihilation operator $A_{i}(t)$ represents a single channel of quantum noise input. $\mathbf{\Lambda}$ is a scattering operator between channels. Both $\mathbf{A}(t)$ and $\mathbf{A}(t)^\ast$ construct a quantum version of Brownian motion processes, while on the other hand, $\mathbf{\Lambda}$ can be thought as a quantum version of a Poisson process \cite{Parthasarathy1992}. 
\\
In the context of open quantum system dynamics, any system observable $X$ will evolve according to 
\begin{dmath}
	X(t) \equiv  U^{\dagger}(t)\left(X \otimes 1 \right)U(t). \label{eq:UnitaryMap}
\end{dmath}
where $1$ is identity operator on $\Gamma\qty(\mathscr{H})^{\otimes ^n}$. Correspondingly, as the analog of \eqref{eq:ClosedSystemHeisenberg}, for an open quantum system, the corresponding Heisenberg equation of motion for a system operator $X$ is given by \cite{gough2009series},
\begin{dmath}
	dX(t) = \left(-i\left[X(t),\mathbb{H}(t)\right] + \frac{1}{2} \mathbb{L}(t)^{\dagger} \left[X(t),\mathbb{L}(t)^\top\right] + \frac{1}{2} \left[\mathbb{L}(t)^{\dagger},X(t)\right]\mathbb{L}(t)\right)dt + d\mathbf{A}^{\dagger}(t)\mathbb{S}^{\dagger}(t)\left[X(t),\mathbb{L}(t)^\top\right] + \left[\mathbb{L}^{\dagger}(t),X(t)\right]\mathbb{S}(t) d\mathbf{A}(t)\nonumber
	+ \text{tr}\left[\left(\mathbb{S}^{\dagger}(t) X(t) \mathbb{S}(t) - X(t)\right)d\bm{\Lambda}(t)^{\top}\right]. \label{eq:QSDE_X}
\end{dmath}
We call this equation a quantum stochastic differential equation (QSDE) for the system observable $X$. All operators in \eqref{eq:QSDE_X} evolve according to  \eqref{eq:UnitaryMap}; i.e.,  $\mathbb{L}(t) \equiv U^{\dagger}(t)\left(\mathbb{L} \otimes 1 \right)U(t)$. 
\\
For the sake of simplicity, if the context is clear, we will drop the time index $t$  from the system observables in the Heisenberg picture; e.g., $X \equiv X(t)$. 

\subsection{Nonlinear Functions of $q$ and $p$, Derivatives and Zero Integrals}

In this subsection, we will examine the QSDEs for a single mode quantum system of the form \eqref{eq:HeisenBergNonLinearQSDE}. Notably, we consider that both $\mathbf{f}$ and $\mathbf{g}$ belong to a class of functions which can be described as a power series in the pair $q,p$ corresponding to the position and momentum operators of the system respectively. We will also introduce a derivative and integration of this class of nonlinear
functions. The generalization of this class to multiple modes along with its derivative and integration will be presented in Section \ref{sec:MultipleModes}. 
\begin{defn}
	A function $f(q,p)$ belongs to the class of nonlinear functions $\mathcal{P}_{q,p}$  if it can be written as a power series of $q$ and $p$ as follows:  
	\begin{dmath}
		X 
		= \sum_{j=1}^{\infty} a_j \phi_j\qty(q,p),{\;\;\;a_j \in \mathbb{C}, \quad k_j, l_j \in \mathbb{N} } \label{eq:P_qp_polynomialExpansion}
	\end{dmath}
	where $\phi_j\qty(q,p)=q^{k_j}p^{l_j}$.
	Also a function f belongs to $\mathcal{P}_q$ or $\mathcal{P}_p$ if it can be written as a power series of $q$ and $p$ respectively as follows :
	\begin{dmath}
		{X 
			= \sum_{j=1}^{\infty} a_j q^{k_j}, \quad Y 
			= \sum_{j=1}^{\infty} b_j p^{l_j}, \quad a_j,b_j \in \mathbb{C}, k_j,l_j \in \mathbb{N}}.
	\end{dmath}
\end{defn}
We notice that under \eqref{eq:P_qp_polynomialExpansion}, $\mathcal{P}_{q,p}$ is a vector space over the complex number field $\mathbb{C}$ , with a basis given by,
\begin{dmath}
	\Phi = \qty{1, q, p, q^2, p^2, qp, q^3,p^3, qp^2,q^2p,\cdots} \label{eq:PqpBasisStandard}.
\end{dmath}
By using the canonical commutation relation  $\comm{q}{p} = i1$, we can always write any $X\in \mathcal{P}_{q,p}$ in `q-p' order as \eqref{eq:P_qp_polynomialExpansion}. Therefore, since both $1$ and $qp$ are already members of the basis $\Phi$,  there is no need to include $pq$ as a member of the basis $\Phi$. Higher order functions of $q$ and $p$ can also be obtained as a linear functions of the members of the basis $\Phi$.
\\
Clearly $\mathcal{P}_q$ and $\mathcal{P}_p$ are subspaces of $\mathcal{P}_{q,p}$ spanned by $\Phi_q = \qty{1, q,  q^2, q^3,\cdots}$ and $\Phi_p = \qty{1, p, p^2, p^3, \cdots}$, respectively.
Suppose we select a set $\Phi_{\mathcal{A}} \subseteq \Phi$. The projection of $X \in \mathcal{P}_{q,p}$ onto the space $\mathcal{A}$ spanned by $\Phi_{\mathcal{A}}$ is denoted by $ \projection{X}{\mathcal{A}}$ . 
The projection of $X \in \mathcal{P}_{q,p}$ onto $\mathcal{A}$, is then given by
\begin{dmath}
	\projection{X}{{\mathcal{A}}} =  \sum_{j=1}^{\infty} a_j' \phi_j\qty(q,p),
\end{dmath}
where $a_j' = a_j$ if $\phi_j\qty(q,p) \in \Phi_{\mathcal{A}}$, and $a_j' = 0$ if $\phi_j\qty(q,p) \notin \Phi_{\mathcal{A}}$.
As an example, let $ \Phi_{\mathcal{A}} = \Phi_q$. The projection of $X$ of the form \eqref{eq:P_qp_polynomialExpansion} onto $\mathcal{P}_q$, is given by
$	\projection{X}{\mathcal{P}_q} =  \sum_{j=1}^{\infty} a_j' \phi_j\qty(q,p)$ ,
where $a_j' = a_j$ if $\phi_j\qty(q,p) \in \Phi_q$, and $a_j' = 0$ if $\phi_j\qty(q,p) \notin \Phi_q$. 
\\
Furthermore, the subspace $\quotient{\mathcal{P}_{q,p}}{\mathcal{A}}$ is defined to be the subspace spanned by $\Phi \setminus {\Phi_{\mathcal{A}}} = \Phi \cap {\Phi_{\mathcal{A}}}^\complement =  {\Phi_{\mathcal{A}}}^\complement$; i.e., the elements of $\Phi$ not in $\Phi_{\mathcal{A}}$, where the complement is taken with $\Phi$ as the universe.

\begin{defn}
The self adjoint subset of $\mathcal{P}_{q,p}$, denoted by $\mathcal{O}_{q,p}$, is set of all $X(q,p) \in \mathcal{P}_{q,p}$ such that $X(q,p) \equiv X(q,p)^\ast$. 
\end{defn}
Using \eqref{eq:P_qp_polynomialExpansion}, we can write any $X \in \mathcal{O}_{q,p}$ as
\begin{dmath}
	X = \sum_{j=1}^{\infty} \frac{1}{2} \qty[a_j \phi_j\qty(q,p) + a_j^\ast \phi_j\qty(q,p)^\ast]  = \sum_{j=1}^{\infty} \real{\qty(a_j)} \frac{\phi_j\qty(q,p)+\phi_j\qty(q,p)^\ast}{2} + \imaginary{\qty(a_j)} i\frac{\phi_j\qty(q,p)^\ast-\phi_j\qty(q,p)}{2} = \sum_{k=1}^{\infty} b_k \varphi_k(q,p), \; {b_k \in \mathbb{R}, \varphi_k(q,p)=\varphi_k(q,p)^\ast}.  \label{eq:O_qp_polynomialExpansion}
\end{dmath}
Observe that $\mathcal{O}_{q,p}$ is a vector space over the field $\mathbb{R}$ , where the collection  $\varphi_k\qty(q,p)$ define a standard basis for the space of functions $\mathcal{O}_{q,p}$.
\\
Now we define derivatives on the space $\mathcal{P}_{q,p}$; see also [\textsection 3]\cite{ARNOLD2008}
\begin{defn}\label{def:DifferentiationOnP}
	For any $X \in \mathcal{P}_{q,p}$, define
	\begin{dgroup}
		\begin{dmath}
			\pdv{X}{p}
			\triangleq \frac{1}{i} \commutator{q}{X}
		\end{dmath}
		\begin{dsuspend}
			and
		\end{dsuspend}
		\begin{dmath}
			\pdv{X}{q} \triangleq \frac{-1}{i} \commutator{p}{X}.
		\end{dmath} \label{eqs:DifferentiationOnP}
	\end{dgroup}
\end{defn}
Both  $\pdv*{q}$ and $\pdv*{p}$  define a surjective  mapping from $\mathcal{P}_{q,p}$ to $\mathcal{P}_{q,p}$, where $\ker\qty(\pdv*{q})  = \mathcal{P}_p$ and $\ker \qty(\pdv*{p}) = \mathcal{P}_q$. 
We will use the following lemma to obtain a general expression for the derivatives defined in this definition.
\begin{lem}\label{lem:DerivationQandP}
	For any $m \geq 1, m \in \mathbb{N}$, we have $\pdv*{q^m}{q} = m q^{m-1}$, and $\pdv*{p^m}{p} = m p^{m-1}$.
\end{lem}

Using the above lemma, for $X$ defined in \eqref{eq:P_qp_polynomialExpansion},  we have
\begin{dgroup}
	\begin{dmath}
		\pdv{X}{q} = \sum_{j=1,   \\ \{j:k_j\geq 1\} }^{\infty} a_j k_j q^{k_j -1}p^{l_j}, 
	\end{dmath}
	\begin{dmath}
		\pdv{X}{p} = \sum_{j=1, \\ \{j:l_j\geq 1\} }^{\infty} a_j l_j q^{k_j}p^{l_j-1}. 
	\end{dmath}
\end{dgroup}
We also notice that the derivatives defined in Definition \ref{def:DifferentiationOnP} preserve self adjointness. That is, if $f \in \mathcal{O}_{q,p}$, then $\pdv*{f}{q}$ and $\pdv*{f}{p}$ belong to $\mathcal{O}_{q,p}$. In the following propositions, we will present some additional properties of the derivatives defined in this definition.
\begin{prop}\label{prp:BasicPropetiesDerivationOnP}
	The derivatives $\pdv*{p}$ and $\pdv*{q}$  have the following properties :
	\begin{enumerate}[label={\alph*)}]
		\item The derivatives $\pdv*{p}$ and $\pdv*{q}$ are commutative operations on $\mathcal{P}_{q,p}$; i.e. ,
		\begin{dmath} 
			\pdv{f}{p}{q} = \pdv{f}{q}{p} {, \; \forall f \in \mathcal{P}_{q,p}}
		\end{dmath}.\label{prt:CommutativeRule}
		\item The derivatives $\pdv*{p}$ and $\pdv*{q}$ satisfy the product rule; i.e., ${\forall X,Y \in \mathcal{P}_{q,p}}$
		\begin{dgroup}
		\begin{dmath} 
			{\pdv{XY}{q} = X \pdv{Y}{q} + \pdv{X}{q}Y,}
		\end{dmath}
		\begin{dsuspend}
			and
		\end{dsuspend}
		\begin{dmath} 
			{\pdv{XY}{p} = X \pdv{Y}{p} + \pdv{X}{p}Y}.
		\end{dmath}
		\end{dgroup}.\label{prt:ProductRule}
		\item The derivatives $\pdv*{p}$ and $\pdv*{q}$ are symmetric; i.e.,
		\begin{dgroup}
		\begin{dmath} 
			{\qty(\pdv{X}{q})^\ast = \pdv{X^\ast}{q},}
		\end{dmath}
		\begin{dsuspend}
			and
		\end{dsuspend}
		\begin{dmath} 			
			{\qty(\pdv{X}{p})^\ast = \pdv{X^\ast}{p}.}
		\end{dmath}
		\end{dgroup}
\label{prt:SymmetricDerivations}
	\end{enumerate}
\end{prop}

In the context of C$^\ast$- algebras, functions on an algebra satisfying the product rule and the symmetric property as in Proposition \ref{prp:BasicPropetiesDerivationOnP}\ref{prt:ProductRule} and \ref{prt:SymmetricDerivations} are called  symmetric `derivations' of the C$^\ast$- algebra; see \cite[Definition 3.2.21]{bratteli1981operator}.
\\
We now define the anti-derivative of a function which belongs to $\mathcal{P}_{q,p}$. 
\begin{defn}
	A function $F \in \mathcal{P}_{q,p}$ is called an anti-derivative of $X \in \mathcal{P}_{q,p}$ with respect to $q$ or $p$ if its derivative with respect to $q$ or $p$ is given by $X$ respectively.  We denote the anti-derivative with respect to $q$ or $p$ by  $\int \cdot dq$ or $\int \cdot dp$ respectively.
\end{defn}
It is obvious that if $F$ is an anti-derivative of $X$ with respect to  $q$ or $p$, then $F + R(p)$ or $F + \tilde{R}(q)$ are also anti-derivatives of $X$ with respect to  $q$ or $p$, for any $R(p) \in \mathcal{P}_q, \tilde{R}(q) \in \mathcal{P}_q$.
\\
Now we define the operations $\int_o \cdot dq $ and $\int_o \cdot dp $, which we refer to as the zero integrals with respect to $q$ and $p$. 
\begin{defn}\label{def:ZeroIntegralQP}
	The zero integrals $\int_o \cdot dq $ and $\int_o \cdot dp $ are mappings from $\mathcal{P}_{q,p} $ to $\quotient{\mathcal{P}_{q,p}}{\mathcal{P}_{q}} $ and from $\mathcal{P}_{q,p} $ to $\quotient{\mathcal{P}_{q,p}}{\mathcal{P}_{p}} $ respectively, such that $\pdv*{\qty(\int_o X dq)}{q} = X$ and $\pdv*{\qty(\int_o X dp)}{p} = X$.
\end{defn}
According to this definition, the zero integral $\int_o \cdot dq$ maps from $\mathcal{P}_{qp}$ to $\mathcal{P}_{qp}\backslash \mathcal{P}_{q}$.  This definition implies that for any $X \in \mathcal{P}_{qp}$, we can always write $\int X dq = \int_o X dq + R(p)$, where $R(p) \in \mathcal{P}_{p}$. This way, $ \int_o X dq$ is the anti-derivative of $X$ with respect to $q$ that has no component in the kernel of $\pdv*{}{q}$, and that $\pdv*{\qty(\int_o X dq)}{q} = X$.
\\
\begin{lem}\label{lem:zero_int_unique}
	The zero integrals  $\int_o \cdot dq $ and $\int_o \cdot dp $ are unique.
\end{lem}

The quantities $\int_o \cdot dq $ and $\int_o \cdot dp $ can be constructed as follows. Let $X \in \mathcal{P}_{q,p}$. Then with $\phi_j(q,p) = q^{k_j}p^{l_j}$, we can write $X$ in \eqref{eq:P_qp_polynomialExpansion} in the following form
\begin{dmath}
	X = \sum\limits_{j=1}^{\infty} a_j F_j(q)G_j(p)
\end{dmath}
where $F_j(q) = {q}^{k_j}$ and $G_j(p) = {p}^{l_j}$ and $a_j \in \mathbb{C}$. From Lemma \ref{lem:DerivationQandP}, and the uniqueness of $\int_o \cdot dq $ and $\int_o \cdot dp $, it follows that $\int_o F_j(q) dq = \frac{1}{k_j+1} {q}^{k_j+1}$ and $\int_o G_j(p) dp  = \frac{1}{l_j+1} {p}^{l_j+1}$.
It follows from Proposition \ref{prp:BasicPropetiesDerivationOnP}\ref{prt:ProductRule}, that if we define
\begin{dgroup}
	\begin{dmath}
		\int_o X dq = \sum\limits_{j=1}^{\infty} a_j \qty(\int_o F_j(q) dq) G_j(p) = \sum\limits_{j=1}^{\infty} a_j \frac{1}{k_j+1} {q}^{k_j+1} {p}^{l_j}
	\end{dmath}
	\begin{dsuspend}
		and,
	\end{dsuspend}
	\begin{dmath}
		\int_o X dp = \sum\limits_{j=1}^{\infty} a_j F_j(q) \qty(\int_o G_j(p)dp)  = \sum\limits_{j=1}^{\infty} a_j \frac{1}{l_j+1} {q}^{k_j} {p}^{l_j+1},
	\end{dmath}
	\label{eqs:Integration}
\end{dgroup}
then we will obtain 
\begin{dmath*}
	\pdv*{\qty(\int_o X dq)}{q} = X$ and $\pdv*{\qty(\int_o X dp)}{p} = X, 
\end{dmath*}
where $\int_o X dq \in \quotient{\mathcal{P}_{q,p}}{\mathcal{P}_p}$ and $\int_o X dp \in \quotient{\mathcal{P}_{q,p}}{\mathcal{P}_q}$. 
\\
As the derivatives $\pdv*{q}$ and $\pdv*{p}$ are commutative in $\mathcal{P}_{q,p}$, from \eqref{eqs:Integration} it follows directly that $\int_o \int_o \cdot dq dp = \int_o \int_o \cdot dp dq$. It also directly follows from \eqref{eqs:Integration} that the zero integrals are symmetric. That is 
\begin{dmath*} 
	{\qty(\int_o X dq)^\ast = \int_o X^* dqx} \text{, and }\; 	{\qty(\int_o X dp)^\ast = \int_o X^* dp.}
\end{dmath*}
%
%

The following proposition summarizes these properties of the zero integrals,
\begin{prop}\label{prp:BasicPropetiesZeroIntegrationOnP}
The zero integrals defined in Definition \ref{def:ZeroIntegralQP} have the following properties:
	\begin{enumerate}[label={\alph*)}]
		\item The zero integrals  $\int_o \cdot dq $ and $\int_o \cdot dp $ are unique.\label{prt:UniquenessOfZI}
		\item The zero integrals  $\int_o \cdot dq $ and $\int_o \cdot dp $ are commutative on $\mathcal{P}_{q,p}$; i.e.,
		\begin{dmath} 
			\int_o \int_o X dp dq = \int_o \int_o X dq dp.
		\end{dmath}\label{prt:CommutativityRuleOfZI}
		\item The zero integrals  $\int_o \cdot dq $ and $\int_o \cdot dp $ are symmetric; i.e.,
		\begin{dgroup}
		\begin{dmath} 
			\qty(\int_o X dq)^\ast = \int_o X^* dq
		\end{dmath}
		\begin{dsuspend}
			and
		\end{dsuspend}
		\begin{dmath} 			
			\qty(\int_o X dp)^\ast = \int_o X^* dp.
		\end{dmath}
		\end{dgroup}\label{prt:SymmetricZI}
	\end{enumerate}
\end{prop}

\subsection{Functions of a Finite Collections of Observables, Their Derivatives and Zero Integrals}\label{sec:MultipleModes}
We now generalize the definitions and results of the previous subsection to the case of multiple mode quantum systems. Let $\mathbf{x}$ be a vector of system observables which do not necessarily commute with each other. In particular, suppose there are $m$ modes in the system and assume $\mathbf{x} = \left[\mathbf{q}^\top \; \mathbf{p}^\top\right]^\top$, where $\mathbf{q}^\top = \left[q_{1} \cdots  q_{m} \right]^\top$ is the vector of position operators for different modes in the system, and $\mathbf{p}^\top = \left[p_{1} \cdots  p_{m} \right]^\top$ is the vector of momentum operators respectively. From the canonical commutation relations, we have for all $t\geq 0$,  and $k,l \leq m$,
	\begin{dmath}
		{\dfrac{1}{i}\commutator{q_{k}}{q_{l}} = 0 \;, \dfrac{1}{i}\commutator{p_{k}}{p_{l}} = 0, \dfrac{1}{i}\commutator{q_{k}}{p_{l}} = \delta_{k-l} }. \label{eqs:CCR_MultiMode}
	\end{dmath}
\begin{defn}
	A function $f\qty(x_1,\cdots,x_{2m})$ belongs to the class $\mathcal{P}_{\qty(x_1,\cdots,x_{2m})} = \mathcal{P}_{\mathbf{x}}$ if it can be written in the form of the following power series:
	\begin{dmath}
		f = \sum_{j=1}^{\infty} a_j x_1^{k_{1,j}}x_2^{k_{2,j}}\cdots x_{2m}^{k_{2m,j}}, {\;\;\;a_j \in \mathbb{C}, \quad k_{i,j} \in \mathbb{N} }.
		\label{eq:P_x_polynomialExpansion}
	\end{dmath}
\end{defn}
Observe that, no matter what order the observables appear in the definition of a function $f$, we can always  rearrange it in the order $x_1,x_2,  \cdots,  x_{2m}$ by suitable use of the commutative relations \eqref{eqs:CCR_MultiMode}.
In \eqref{eq:P_x_polynomialExpansion}, the collection of polynomial functions
\begin{dmath}
	\phi_j\qty(x_1,x_2,\cdots, x_{2m}) = x_1^{k_{1,j}}x_2^{k_{2,j}}\cdots x_{2m}^{k_{2m,j}},
\end{dmath}
define a standard basis for the space of functions $\mathcal{P}_{\mathbf{x}}$; i.e., we write
\begin{align*}
\Phi  = &\left\{ 1, \right.\\
&x_1, x_2, \cdots, x_{2m},\\
&x_1^2,x_2^2, \cdots, x_{2m}^2,\\
&x_1x_2, x_1 x_3, \cdots, x_1x_{2m},\\
&x_2x_3, x_2 x_4, \cdots, x_2x_{2m}, \\
& \cdots\\
& \vdots\\
& \left. \cdots,\right\}.
\end{align*}
This, in turn makes $\mathcal{P}_{\mathbf{x}}$ a vector space over the complex numbers $\mathbb{C}$.
\\
As for the single mode case, we will use the following definition of the derivatives on $\mathcal{P}_{\mathbf{x}}$
\\
\begin{defn}
	For any $X \in \mathcal{P}_{\mathbf{x}}$, define
	\begin{dgroup}
		\begin{dmath}
			\pdv{X}{p_{i}}
			\triangleq \frac{1}{i} \commutator{q_{i}}{X},
		\end{dmath}
		\begin{dmath}
			\pdv{X}{q_{i}} \triangleq \frac{-1}{i} \commutator{p_{i}}{X}.
		\end{dmath} \label{eqs:DifferentiationOnPMulti}
	\end{dgroup}
\end{defn}
Using a similar approach as in Proposition \ref{prp:BasicPropetiesDerivationOnP}, one can verify that these derivatives are commutative on $\mathcal{P}_{\mathbf{x}}$, symmetric, and satisfy the product rule. 

\begin{defn}
	The self adjoint subset of $\mathcal{P}_{\mathbf{x}}$, denoted as $\mathcal{O}_{\mathbf{x}}$, is the collection of all $X \in \mathcal{P}_{\mathbf{x}}$ such that $X\qty(x_1,\cdots,x_{2m}) \equiv X\qty(x_1,\cdots,x_{2m})^\ast$. 
\end{defn}
In a similar way to \eqref{eq:O_qp_polynomialExpansion}, we can construct a set of functions $\varphi_k\qty(\mathbf{x})$ which define a standard basis for $\mathcal{O}_{\mathbf{x}}$. 
Suppose we select $\Phi_{\mathcal{A}} \subseteq \Phi$. The projection of $f \in \mathcal{P}_{\mathbf{x}}$ onto a space $\mathcal{A}$ spanned by $\Phi_{\mathcal{A}}$, is given by
$	\projection{f}{{\mathcal{A}}} =  \sum_{j=1}^{\infty} a_j' \phi_j\qty(q,p)$,
where $a_j' = a_j$ if $\phi_j \in \Phi_{\mathcal{A}}$, and $a_j' = 0$ if $\phi_j \notin \Phi_{\mathcal{A}}$. Furthermore, the subspace $\quotient{\mathcal{P}_{\mathbf{x}}}{\mathcal{A}}$ is defined to be the subspace spanned by $\Phi \setminus {\Phi_{\mathcal{A}}} = \Phi \cap {\Phi_{\mathcal{A}}}^\complement$.

Suppose $\mathbf{z} = \qty{z_1,\cdots,z_k} \subseteq \qty{x_1,\cdots,x_{2m}} $. The space $\mathcal{P}_{\mathbf{z}}$ is the set of all functions that can be written as a power series of $\qty{z_1,\cdots,z_k}$. We also define $	\bar{\mathcal{P}}_{\mathbf{z}}	\triangleq  \mathcal{P}_{\mathbf{z}^\complement}$.
That is, $\bar{\mathcal{P}}_{\mathbf{z}}$ is the set of functions of $\qty{x_1 \cdots x_{2m} }\setminus \mathbf{z}$ .
Also notice  that the following properties hold:
\begin{dmath*}
	{\bar{\mathcal{P}}_{\{\}} = \mathcal{P}_{\mathbf{x}}, \; \bar{\mathcal{P}}_{\mathbf{x}} = \mathbb{C}}.
\end{dmath*}
Using this notation, we observe that for any $ 1 \leq i \leq 2m$,
\begin{dmath*}
	\quotient{\bar{\mathcal{P}}_{\qty(x_1,\cdots,x_{i-1})}}{\bar{\mathcal{P}}_{\qty(x_1,\cdots,x_{i})}}	\bigcap  \quotient{\bar{\mathcal{P}}_{\qty(x_1,\cdots,x_{i})}}{\bar{\mathcal{P}}_{\qty(x_1,\cdots,x_{i+1})}} = \{0\}.
\end{dmath*}

We then define the zero integral as follows.
\begin{defn}\label{def:ZeroIntegralsGeneral}
	The zero integral $\int_o \cdot dx_i $ is a mapping from $\mathcal{P}_{\mathbf{x}} $ to $\quotient{\mathcal{P}_{\mathbf{x}}}{\bar{\mathcal{P}}_{x_i}} $, such that $\pdv*{\qty(\int_o X dx_i)}{x_i} = X, \forall X \in \mathcal{P}_{\mathbf{x}}$.
\end{defn}
We now observe that from the properties of the derivatives $\pdv*{}{x_i}$, for any two operators $X,X' \in \quotient{\mathcal{P}_{\mathbf{x}}}{\bar{\mathcal{P}}_{x_i}}$, the condition $\pdv*{\qty(X-X')}{dx_i}=0$ will only be satisfied if $X=X'$. Therefore, as in the previous subsection this leads to the uniqueness of $\int_o dx_i$. The commutative and symmetric properties of $\int_o dx_i$ can be concluded directly from the commutativity of $\pdv*{x_i}$ on $\mathcal{P}_{\mathbf{x}}$ and Definition \ref{def:ZeroIntegralsGeneral}.
\\
Using the product rule and $\pdv*{x_i}{x_i} = 1$, in a similar way as in Lemma \ref{lem:DerivationQandP}, we obtain $\pdv*{x_i^m}{x_i} = mx_i^{m-1}$. Therefore, we can readily verify that for $X$ of the form  \eqref{eq:P_x_polynomialExpansion}
\begin{dgroup}
	\begin{dmath}
		\pdv{X}{x_i} = \sum_{j=1,   \\ \{j:k_{i,j}\geq 1\} }^{\infty} a_j k_{i,j} x_1^{k_{1,j}} \cdots x_i^{k_{i,j}-1}\cdots x_{2m}^{k_{2m,j}},
	\end{dmath}
	\begin{dsuspend}
		and
	\end{dsuspend}
	\begin{dmath}
		\int_o X dx_i = \sum_{j=1}^{\infty} a_j x_{1}^{k_{1,j}} \cdots \int_o x_{i}^{k_{i,j}} dx_i \cdots x_{2m,j}^{k_{2m,j}} = \sum_{j=1}^{\infty} a_j \dfrac{1}{k_{i,j}+1} x_{1}^{k_{1,j}} \cdots  x_{i}^{k_{i,j}+1} \cdots x_{2m}^{k_{2m}}.
	\end{dmath}
\end{dgroup}
We can also construct a projection of the zero integral onto a subspace of $\mathcal{P}_{\mathbf{x}}$. For example, the integration term $\left. \int_o \pdv{f}{x_i} d x_i \right|_{\mathcal{P}_{(x_i,\cdots,x_{2m})}}$ is the projection of $\int_o \pdv{f}{x_i} d x_i$ onto $\mathcal{P}_{\qty(x_i,\cdots,x_{2m})}$, which in turn maps from $\mathcal{P}_{\mathbf{x}} $ to $\left.\qty(\quotient{\mathcal{P}_{\mathbf{x}}}{\bar{\mathcal{P}}_{x_i}}) \right|_{\bar{\mathcal{P}}_{x_1,\cdots,x_{i-1}}}  = \quotient{\bar{\mathcal{P}}_{x_1,\cdots,x_{i-1}}}{\bar{\mathcal{P}}_{x_1,\cdots,x_i}}$. 
\\
In the following lemma, we will show that we can expand any function $f \in \mathcal{P}_{\mathbf{x}}$ as a series of zero integrals. We will use this lemma in the next subsection to prove necessary and sufficient conditions for a vector whose elements belong to $\mathcal{P}_{\mathbf{x}}$ to be the gradient of a potential. Before that, we recall that a permutation $\sigma$ of the set $\qty{1,\cdots,2m}$ is a bijective mapping from $\qty{1,\cdots,2m}$ onto itself. 
\begin{lem} \label{lem:Py_Expansion}
	Any $f \in \mathcal{P}_{\mathbf{x}}$ can be expanded as a series of integrals with respect to a permutation of $\mathbf{x}$, $\mathbf{y} = \left[x_{\sigma\qty(1)} \cdots x_{\sigma\qty(2m)}\right]$, 
	\begin{dmath}
		f = \int_o \pdv{f}{y_1} d y_1 + \left. \int_o \pdv{f}{y_2} d y_2 \right|_{\mathcal{P}_{\qty(y_2,\cdots,y_{2m})}} + \cdots + \left. \int_o\pdv{f}{y_i} d y_i \right|_{\mathcal{P}_{\qty(y_i,\cdots,y_{2m})}}+ \cdots + \left. \int_o\pdv{f}{y_{2m}} d y_{2m} \right|_{\mathcal{P}_{y_{2m}}} + C = \int_o \pdv{f}{y_1} d y_1 + R_{y_1}, \label{eq:f_expandY}
	\end{dmath}
	where $R_{y_1} \in \bar{\mathcal{P}}_{y_1}$, and $C \in \mathbb{C}$ is a constant.
\end{lem}

\subsection{The Gradient Vector on $\mathcal{P}_{\mathbf{x}}$}
In this subsection, we will introduce the concept of a gradient vector on $\mathcal{P}_{\mathbf{x}}$. Let $\mathcal{P}_{\mathbf{x}}^{k \times l}$ be a space of $k\times l$ matrices whose elements belong to $\mathcal{P}_{\mathbf{x}}$.
\begin{defn}
	A vector $\mathbf{g} \in \mathcal{P}_{\mathbf{x}}^{2m \times 1}$ is a gradient of $f$ with respect to $\mathbf{x}$, if $\pdv*{f}{\mathbf{x}} = \mathbf{g}^\top$. Furthermore, $\mathbf{g}$ is a \emph{gradient} with respect to $\mathbf{x}$ if there exists an $f \in \mathcal{P}_{\mathbf{x}}$ such that $\pdv*{f}{\mathbf{x}} = \mathbf{g}^\top$.
\end{defn}
For the sake of simplicity, we also use the following notation for the integral expansion in \eqref{eq:f_expandY},
\begin{dmath}
	f = \int_o \pdv{f}{\mathbf{y}} d \mathbf{y} + C . 
\end{dmath}
In vector calculus, a gradient defined as above is also known as a conservative vector field. In simply connected spaces like $\mathbb{R}^2$ and $\mathbb{R}^3$, a gradient can be characterized by the fact that the curl operation acting on it will equal zero \cite{Stewart2011}. 
\\
The following theorem shows that a gradient with respect to $\mathbf{x}$ must be expandable by a series of zero integrals of permutations of $\mathbf{x}$. 
\begin{thm}\label{thm:gradient}
	Consider $\mathbf{g} \in \mathcal{P}_{\mathbf{x}}^{2m \times 1}$ 	and $f \in \mathcal{P}_{\mathbf{x}}$. Then $\mathbf{g}$ is a gradient of $f$ with respect to $\mathbf{x}$ if and only if $f$ can be expanded in the integral sum form \eqref{eq:f_expandY}, using elements of $\mathbf{g}$ for any permutation of $\mathbf{x}$.
\end{thm}

In two-dimensional vector calculus, Green's theorem leads to the fact that every curl zero vector (irrotational vector) is a conservative vector field in a simply connected region \cite{Stewart2011}. It is interesting to consider whether every curl zero vector in $\mathcal{P}_{\mathbf{x}}$ is also a gradient. The following theorem shows that the answer to this question is affirmative.
\begin{thm}\label{thm:Curl}
	Consider $\mathbf{g} \in \mathcal{P}_{\mathbf{x}}^{2m \times 1}$. 
	Then $\mathbf{g}$ is a gradient with respect to $\mathbf{x}$ if only if, 
	\begin{dmath}
		\pdv{g_i}{x_j} = \pdv{g_j}{x_i} {, \forall i \neq j, 1 \leq i,j \leq 2m}. \label{eq:CurlGeneralized}
	\end{dmath}
\end{thm}

The following corollary is a straightforward generalization of Theorem \ref{thm:gradient} to a set of observables.
\begin{cor}
	Let $\mathbf{g} \in \mathcal{P}_{\mathbf{x}}^{2m \times k}$ 	and $\mathbf{f} \in \mathcal{P}_{\mathbf{x}}^{k \times 1}$ be given. Then $\mathbf{g}$ is a gradient of $\mathbf{f}$ with respect to $\mathbf{x}$ if and only if each element of $\mathbf{f}_i$ can be expanded using elements of $\mathbf{g}$ for any permutation of $\mathbf{x}$.	
\end{cor}

\section{Main Results}\label{sec:MainResults}
In this section, we will give sufficient and necessary conditions for a QSDE corresponding to a single mode quantum system interacting with single environment field to be physically realizable. For this purpose, let $\mathbf{x} = \left[q \; p\right]^\top$.
\begin{defn}\label{def:PhysRealizableQSDESingleMode}
	The QSDE in \eqref{eq:HeisenBergNonLinearQSDE} is said to be physically realizable if there exists a pair $\mathbb{H} \in \mathcal{O}_{q,p}$, and  $\mathbb{L} \in \mathcal{P}_{q,p}$ satisfying the following equations
	\begin{dgroup}\label{eqs:NonlinearQSDESingleMode}
		\begin{dmath}
			\mathbf{f}(\mathbf{x}) =  \mathcal{L}(\mathbf{x}) = -i\left[\mathbf{x},\mathbb{H}\right] + \frac{1}{2}\left[\qty(\mathbf{g}(\mathbf{x})^\ast\mathbb{L})^\ast +  	\mathbf{g}(\mathbf{x})^\ast\mathbb{L}\right], \label{eq:quantumMarkovGenerator}
		\end{dmath}
		\begin{dmath}
			\mathbf{g}(\mathbf{x}) =   \left[\mathbf{x},\mathbb{L}\right]. \label{eq:G}
		\end{dmath}
	\end{dgroup}	
\end{defn}
As Definition \ref{def:PhysRealizableQSDESingleMode} implies, the QSDE in \eqref{eq:HeisenBergNonLinearQSDE} is physically realizable if it is corresponds to dynamic of an open quantum system with a coupling operator $\mathbb{L}$ and a Hamiltonian $\mathbb{H}$ where the Hamiltonian $\mathbb{H}$ needs to be self-adjoint. 


Before we state our main result, we will introduce a notion of a commutator-conservative mapping $\mathbf{f}$. First, let $\bm{\Sigma} = \begin{bmatrix}
0 && 1 \\ -1 && 0
\end{bmatrix}$.

\begin{defn}
	A mapping  $\mathbf{j} \in \mathcal{P}_{q,p}^{2 \times 1}$ is said to be commutator-conservative 
	if there exists $\mathbb{J} \in \mathcal{P}_{q,p}$, such that  $\mathbf{j} =  -i\left[\mathbf{x},\mathbb{J}\right]$.
\end{defn}
If we consider a closed quantum system with dynamics as in \eqref{eq:ClosedSystemHeisenberg}, the function $\mathbf{f} = \pdv*{\mathbf{x}}{t}$ is commutator-conservative if there exists a Hamiltonian $\mathbb{H}$, such that  $\mathbf{f} =  -i\left[\mathbf{x},\mathbb{H}\right]$. Furthermore, if the QSDE \eqref{eq:HeisenBergNonLinearQSDE} is physically realizable, the function $\mathbf{g}$ will be commutator-conservative for $\mathbb{J} = -i\mathbb{L}$. In what follows, by using results from the preliminary section, we will arrive at a conclusion that a function $\mathbf{f}$ with self-adjoint elements is commutator-conservative if there exists a potential observable $\mathbb{H}$, such that  $\mathbf{f} =  \pdv*{\mathbb{H}}{ \qty(\bm{\Sigma} \mathbf{x})}$, or if it is satisfies $-\pdv{f_1}{q} = \pdv{f_2}{p}$, where $f_1$, and $f_2$ are self-adjoint elements of $\mathbf{f}$. In this way, we could think of a commutator-conservative mapping $\mathbf{f} \in \mathcal{O}_{q,p}^{2\times 1}$ as a conservative vector field of a potential $\mathbb{H}$  with the axis being  $x = p$, and $y = -q$.
\\
The following result gives a sufficient and necessary condition for a mapping $\mathbf{j} \in \mathcal{P}_{q,p}^{2 \times 1}$ to be commutator-conservative. 
\begin{thm}\label{thm:ClosedPhysicalRealizability}
	Consider  $\mathbf{j} = \qty[j_1 \; j_2]^\top \in \mathcal{P}_{q,p}^{2 \times 1}$. $\mathbf{j}$ is commutator-conservative if and only if 
	\begin{dmath}
		\left. \int_o j_1 dp \right|_{ \mathcal{P}_{p}} - \int_o j_2 dq  = \int_o j_1 dp - \left. \int_o j_2 dq \right|_{ \mathcal{P}_{q}}.   \label{eq:ClosedPhysicalRealizability}
	\end{dmath}
\end{thm}
\begin{pf}
	Suppose $\mathbf{j}$ satisfies \eqref{eq:ClosedPhysicalRealizability}. Then we can construct an operator $\mathbb{J} \in \mathcal{P}_{q,p}$, as follows:
	\begin{dmath}
		\mathbb{J} = \left. \int_o j_1 dp \right|_{ \mathcal{P}_{p}} - \int_o j_2 dq + C \label{eq:J_PthenQ_expansion}
	\end{dmath}
	where $C$ is a constant self adjoint operator. 	By 	\eqref{eq:ClosedPhysicalRealizability}, we have $j_1 = \pdv*{\mathbb{J}}{p}$, $j_2 = - \pdv*{\mathbb{J}}{q}$. Therefore,	
	\begin{dmath*}
		{\begin{bmatrix}	j_1 \\ j_2 \end{bmatrix} 
			= \begin{bmatrix}0 && 1 \\ -1 && 0\end{bmatrix} 
			\begin{bmatrix}
				\pdv*{\mathbb{J}}{q} \\ \pdv*{\mathbb{J}}{p} 
			\end{bmatrix} 
			= \frac{1}{i} \begin{bmatrix}
				i \pdv*{\mathbb{J}}{p} \\ -i \pdv*{\mathbb{J}}{q}
		\end{bmatrix} }\\
	{= \frac{1}{i} \begin{bmatrix}
			\left[q,\mathbb{J}\right] \\ 
			\left[p,\mathbb{J}\right] 
		\end{bmatrix} = -i\left[\mathbf{x},\mathbb{J}\right]} . \label{eq:ProvingCommutatorConservativeness}
\end{dmath*}
Now suppose $\mathbf{j}$ is commutator-conservative; i.e, $\mathbf{j} = -i\comm{\mathbf{x}}{\mathbb{J}}$, for some $\mathbb{J} \in \mathcal{P}_{q,p}$. Also, suppose $\mathbb{J}_q = \left. \mathbb{J}\right|_{ \mathcal{P}_{q}}$, and $\mathbb{J}_p = \left. \mathbb{J}\right|_{ \mathcal{P}_{p}}$. Then we can expand $\mathbb{J}$ as follows:
\begin{dmath}
	{\mathbb{J} = \int_o \pdv{ \mathbb{J}}{p} dp + \mathbb{J}_q = \mathbb{J}_p + \int_o \pdv{\mathbb{J}}{q} dq }. \label{eq:J_expansion}
\end{dmath}
From \eqref{eq:ProvingCommutatorConservativeness}, we have  $j_1 = \pdv*{\mathbb{J}}{p}$, $j_2 = - \pdv*{\mathbb{J}}{q}$, and therefore by \eqref{eq:J_expansion}, we have
\begin{dmath*}
\mathbb{J}_q = {\left. \mathbb{J}_q\right|_{ \mathcal{P}_{q}} = \left. \qty(\mathbb{J}_p - \int_o j_2 dq - \int_o j_1 dp )\right|_{ \mathcal{P}_{q}}} = - \left. \int_o j_2 dq \right|_{ \mathcal{P}_{q}} + C
\end{dmath*}
and
\begin{dmath*}
	\mathbb{J}_p = {\left. \mathbb{J}_p\right|_{ \mathcal{P}_{p}} = \left. \qty(\mathbb{J}_q + \int_o j_2 dq + \int_o j_1 dp )\right|_{ \mathcal{P}_{p}}} = \left. \int_o j_1 dp \right|_{ \mathcal{P}_{p}} + C.
\end{dmath*}
Hence, \eqref{eq:ClosedPhysicalRealizability} is satisfied.\qed\end{pf}
\begin{cor}\label{cor:ThreeConditionOfCommutator-Conservativity}
	Consider  $\mathbf{j} \in \mathcal{P}_{q,p}^{2 \times 1}$. The following three conditions are equivalent,
	\begin{enumerate}
		\item $\mathbf{j}$ is commutator-conservative,
		\item $\mathbf{j}$ is a gradient with respect to $\mathbf{y} = \bm{\Sigma} \mathbf{x}$,
		\item $\mathbf{j}$ satisfies,
		\begin{dmath}
			-\pdv{j_1}{q} = \pdv{j_2}{p}. \label{eq:Curl_pq}
		\end{dmath}
	\end{enumerate}
	
\end{cor}

\begin{pf}
	This result is combination of Theorem \ref{thm:ClosedPhysicalRealizability}, Theorem \ref{thm:gradient}, and Theorem \ref{thm:Curl}.
\qed\end{pf}

\begin{prop}\label{prp:SelfAdjointNessOfJ}
	Using the notation of Theorem \ref{thm:ClosedPhysicalRealizability}, if a mapping $\mathbf{j}$ is commutator-conservative and  has self-adjoint elements; i.e., $j_1,j_2 \in \mathcal{O}_{q,p}$, then there exists a self-adjoint operator $\mathbb{J} \in \mathcal{P}_{q,p}$ such that $\mathbf{j} = -i\comm{\mathbf{x}}{\mathbb{J}}$. 
\end{prop}
\begin{pf}
	To see this, without loss of generality assume that $j_1,j_2 \neq 0$. Otherwise we can select $\mathbb{J}$ as a real constant $C \in \mathbb{R}$.  Suppose that there is no self-adjoint operator $\mathbb{J}$ satisfying 
	\begin{subequations}
		\begin{align}
		j_1 = -i\comm{q}{\mathbb{J}}, \\
		j_2 = -i\comm{p}{\mathbb{J}}.
		\end{align}\label{eqs:j1_j2}
	\end{subequations}
	Observe that we can write $\mathbb{J} = \mathbb{J}_r + i \mathbb{J}_i$, where both $\mathbb{J}_r$ and $\mathbb{J}_i$ are self-adjoint operators, and $\mathbb{J}_i \neq 0$.
	Now observe that since $j_1,j_2$ are self-adjoint operators
	\begin{subequations}
		\begin{align}
		0 =& j_1^\ast - j_1 = \qty(-i\qty[q\mathbb{J} - \mathbb{J} q])^\ast - \qty(-i\qty[q\mathbb{J} - \mathbb{J} q]) = -2[q,\mathbb{J}_i], \label{eq:Ji_Pq}\\
		0 =& j_2^\ast - j_2 = \qty(-i\qty[p\mathbb{J} - \mathbb{J} p])^\ast - \qty(-i\qty[p\mathbb{J} - \mathbb{J} p]) = -2[p,\mathbb{J}_i].\label{eq:Ji_Pp}
		\end{align}
	\end{subequations}
	Therefore, from \eqref{eq:Ji_Pq}, $\mathbb{J}_i \in \mathcal{P}_q$, while at the same time from \eqref{eq:Ji_Pp}, $\mathbb{J}_i \in \mathcal{P}_p$. Hence $\mathbb{J}_i \in \mathcal{P}_q \bigcap \mathcal{P}_p = \mathbb{C}$. But if this is the case, then taking $\mathbb{J}' = \mathbb{J}_r$ will give $j_1 = -i\comm{q}{\mathbb{J}'}, j_2 = -i\comm{p}{\mathbb{J}'}$. Hence we have arrived at a contradiction.
\end{pf}
Due to the nature of the unitary operator $U(t)$, for a physically realizable quantum system \eqref{eq:HeisenBergNonLinearQSDE}, if a pair of observables at the initial time has a commutator equal to $\commutator{X(0)}{Y(0)} = i \epsilon, \quad \epsilon \in \mathbb{R}$, then this commutator will remain constant in the future. Mathematically, this is equivalent to the following,
\begin{dmath}
	\commutator{X(t)}{Y(t)} = X(t) Y(t) - Y(t) X(t) = U(t)^\ast X(0) U(t) U(t)^\ast Y(0) U(t) - U(t)^\ast Y(0) U(t) U(t)^\ast X(0) U(t) = U(t)^\ast  \commutator{X(0)}{Y(0)} U(t) = i \epsilon. \label{eq:CommutationPreservation}
\end{dmath}
Due to this property, we say that the commutation relations are preserved.
In particular, for $\mathbf{x}(t) = \left[q(t) \; p(t)\right]^\top$, we have for any $t \geq 0$, $\commutator{\mathbf{x}(t)}{\mathbf{x}(t)^\top} =  \commutator{\mathbf{x}(0)}{\mathbf{x}(0)^\top} = i \bm{\Sigma}$.
For a closed quantum system, if the commutation relations are preserved, then
\begin{dmath}
	{\dfrac{d \left[\mathbf{x}(t),\mathbf{x}(t)^\top\right]}{dt} =  i \dfrac{d \bm{\Sigma} }{dt} = 0_{2 \times 2}}.
\end{dmath}
The following proposition and its corollary show that for any commutator-conservative mapping $\mathbf{f}$, the corresponding closed quantum system \eqref{eq:ClosedSystemHeisenberg} will preserve the commutation relations.

\begin{prop}\label{prp:CommutativityPreservationClosed}
	Consider  $\mathbf{j} \in \mathcal{P}_{q,p}^{2 \times 1}$. If  $\mathbf{j}$ is a commutator-conservative mapping 
	then
	\begin{dmath}
		{\qty[\mathbf{j},\mathbf{x}^\top]+\qty[\mathbf{x},\mathbf{j}^\top] = 0_{2 \times 2}}. \label{eq:CommutatorPreservationJ}
	\end{dmath}
\end{prop}
\begin{pf}	
	Let $\mathbf{j}$ be commutator-conservative. Then by Proposition \ref{prp:BasicPropetiesDerivationOnP}\ref{prt:CommutativeRule} and the existence of $\mathbb{J} \in \mathcal{P}_{q,p}$ by Theorem \ref{thm:ClosedPhysicalRealizability}, we have
	\begin{dmath*}
		{\left[j_1,p\right] = i \pdv{ j_1}{ q} = i \pdv{ \mathbb{J}}{ q}{p} = - i \pdv{p} \qty(-\pdv{\mathbb{J}}{q}) } =  -\left[q,j_2\right].
		\end{dmath*}
		The condition $\left[j_2,q\right] = -\left[q,j_1\right]$ can also be established in the similar way.
	Then, we obtain $\left[j_1,p\right] = -\left[q,j_2\right]$ and $\left[j_2,q\right] = -\left[q,j_1\right]$. However, this is equivalent to \eqref{eq:CommutatorPreservationJ} and hence the proof is completed.
\qed\end{pf}
\begin{cor}
		Consider  $\mathbf{f} \in \mathcal{O}_{q,p}^{2 \times 1}$. If $\mathbf{f}$ is a commutator-conservative mapping,
		 then the corresponding closed quantum system dynamics \eqref{eq:ClosedSystemHeisenberg} preserve the commutation relations.
\end{cor}
\begin{pf}
		First observe that for the system \eqref{eq:ClosedSystemHeisenberg} with  $\mathbf{x} = \left[q \; p\right]^\top$, the preservation of the commutation relation will be satisfied if
$				\pdv*{\comm{\mathbf{x}}{\mathbf{x}^\top}}{dt} = \left[\mathbf{f},\mathbf{x}^\top\right]+\left[\mathbf{x},\mathbf{f}^\top\right] = 0_{2 \times 2} $.
		The result then follows by Proposition \ref{prp:CommutativityPreservationClosed}.
\qed\end{pf}
We will use the following Lemma in the main theorem to show that the constant contribution obtained during integration of $\mathbf{g}$ in \eqref{eq:HeisenBergNonLinearQSDE} to construct the coupling operator $\mathbb{L}$ does not affect the physical realizability property of QSDEs.
\begin{lem}\label{lem:roleOfConstant}
	 Consider the QSDEs given in \eqref{eq:HeisenBergNonLinearQSDE}. Let
	\begin{dgroup}
		\begin{dmath}
			Z = \left. \int_o g_1 dp \right|_{ \mathcal{P}_{p}} - \int_o g_2 dq , \label{eq:Z}
		\end{dmath}
		\begin{dmath}
			\mathbb{L} = -i \qty(	Z + C)  , \label{eq:L}
		\end{dmath}
		\begin{dmath}
			\mathbf{f}_L = \frac{\qty( \mathbf{g}^\ast \mathbb{L})+\qty( \mathbf{g}^\ast \mathbb{L})^\ast}{2}, \label{eq:fL}
		\end{dmath}				
	\end{dgroup}
	where $C \in \mathbb{C}$ is a constant. Suppose $\mathbf{g}$ is a commutator-conservative mapping. Let, $\mathbb{L}_1 = -i \qty(	Z + C_1 )$ and $\mathbf{f}_{L,1} = \frac{1}{2}\qty[\qty( \mathbf{g}^\ast \mathbb{L}_1)+\qty( \mathbf{g}^\ast \mathbb{L}_1)^\ast]$, where $C_1$ is a constant. Then for any $C_1 \in \mathbb{C}$, $\mathbf{f} - \mathbf{f}_{L,1}$ is commutator-conservative if $\mathbf{f} - \mathbf{f}_{L}$ is commutator-conservative.
\end{lem}
\begin{pf}
Suppose	$\mathbf{f} - \mathbf{f}_L$ is commutator-conservative, and let $\mathbf{j} = \mathbf{f}-\mathbf{f}_L$. Then $\mathbf{f}-\mathbf{f}_{L,1} = \mathbf{j} - \mathbf{f}_1$, with $\mathbf{f}_1 = \frac{i}{2}\qty[\mathbf{g} \qty(C_1 - C)^\ast -\mathbf{g}^\ast \qty(C_1 - C)]$ . However $\mathbf{g}$ is commutator-conservative by assumption, which implies the commutator-conservativeness of $\mathbf{f}_1$, hence so is $\mathbf{f}-\mathbf{f}_{L,1}$.	
\qed
\end{pf}

Now using the concept of a commutator-conservative mapping, we will establish a condition for the physical realizability of the open quantum system given in \eqref{eq:HeisenBergNonLinearQSDE}, which is given in the following theorem.
\begin{thm}\label{thm:PhysicalRealizabilityOpen}
	Consider the QSDEs given in \eqref{eq:HeisenBergNonLinearQSDE}. Using the notation of Lemma \ref{lem:roleOfConstant},  the following conditions are equivalent
	\begin{enumerate}
		\item The QSDEs in \eqref{eq:HeisenBergNonLinearQSDE} are physically realizable,
		\item Both $\mathbf{g}$ and $\mathbf{f} - \mathbf{f}_L$ are commutator-conservative mappings for some $C \in \mathbb{C}$ in \eqref{eq:L}, and the elements of $\mathbf{f}$ are self-adjoint, 
		\item Both $\mathbf{g}$ and $\mathbf{f} - \mathbf{f}_L$ are commutator-conservative mappings for any $C \in \mathbb{C}$ in \eqref{eq:L}, and the elements of $\mathbf{f}$ are self-adjoint. 
	\end{enumerate}
\end{thm}
\begin{pf}
	First, assume that the QSDEs \eqref{eq:HeisenBergNonLinearQSDE} are physically realizable. From \eqref{eq:G} with  $\mathbb{J} = -i\mathbb{L}$, we conclude the commutator-conservativeness of $\mathbf{g}$. Therefore, there exists $C \in \mathbb{C}$ such that $\mathbb{L}$ satisfies \eqref{eq:L}. We observe that elements of the function $\mathbf{f}$ in \eqref{eq:quantumMarkovGenerator} are self-adjoint. Further, from $\mathbf{f}_L$ given in \eqref{eq:fL}, we have $\mathbf{f} - \mathbf{f}_L = -i\comm{\mathbf{x}}{\mathbb{H}}$ which is also commutator-conservative. Therefore, this implies the condition (2).
	\\
	Now suppose condition (2) holds. Lemma \ref{lem:roleOfConstant} shows that condition (2) implies condition (3).
	\\
	Lastly, assume that condition (3) of the theorem is satisfied. Using the commutator-conservativeness of $\mathbf{g}$, Theorem \ref{thm:ClosedPhysicalRealizability} implies that there exists an $\mathbb{L} \in \mathcal{P}_{q,p}$ given in \eqref{eq:L}, such that $\mathbf{g} = \comm{\mathbf{x}}{\mathbb{L}}$. Furthermore, we can construct $\mathbf{f}_L$ as in \eqref{eq:fL}, which is a vector with self-adjoint operator elements. 
	Since $\mathbf{f} - \mathbf{f}_L$ is commutator-conservative and self-adjoint for any $C \in \mathbb{C}$ in \eqref{eq:L}, then by Theorem \ref{thm:ClosedPhysicalRealizability} and Proposition \ref{prp:SelfAdjointNessOfJ}, there exists an $\mathbb{H} \in \mathcal{O}_{q,p}$, such that $\mathbf{f} - \mathbf{f}_L = -i\qty[\mathbf{x},\mathbb{H}]$. From the existence of the Hamiltonian $\mathbb{H}$ and the coupling operator $\mathbb{L}$, it follows that both $\mathbf{f}$ and $\mathbf{g}$ satisfy \eqref{eq:quantumMarkovGenerator} and \eqref{eq:G} respectively, and therefore the QSDEs are physically realizable.\qed\end{pf}

From the properties of the open quantum system dynamics \eqref{eq:CommutationPreservation}, QSDEs which correspond to an open quantum system will always preserve the commutation relations. In the following proposition, we will show that  physically realizable QSDEs do indeed preserve the commutation relations. However, before we proceed, we observe that for open quantum dynamics as in \eqref{eqs:NonlinearQSDESingleMode} and for $\mathbf{x} = \left[q \; p\right]^\top$, when quantum randomness effects are taken into account, the  preservation of the commutation relations corresponds to $d \left[\mathbf{x},\mathbf{x}^\top\right] = \mathbf{0}_{2\times 2}$. 
\begin{prop}\label{prp:PreservationCommutativity}
	Consider QSDEs given as in \eqref{eq:HeisenBergNonLinearQSDE}. If these QSDEs are physically realizable, then they preserve the commutation relations.
\end{prop}
\begin{pf}
	Let $\mathbf{f}_C = \mathbf{f} - \mathbf{f}_L$. Using the quantum It\^{o} rule \cite{Parthasarathy1992}, the preservation of the commutation relations for the QSDEs \eqref{eq:HeisenBergNonLinearQSDE} is equivalent to
	\begin{dmath*}
		\mathbf{0}_{2\times 2} = d \left[\mathbf{x},\mathbf{x}^\top\right] =   \qty(\underbrace{\left[\mathbf{f}_C,\mathbf{x}^\top\right]+\left[\mathbf{x},\mathbf{f}_C^\top\right]}_{A_1} + \underbrace{\left[\mathbf{f}_L,\mathbf{x}^\top\right]+\left[\mathbf{x},\mathbf{f}_L^\top\right]}_{A_2} + \underbrace{\left[\mathbf{g}^\ast,\mathbf{g}^\top\right]}_{A_3})dt +
		\underbrace{\left[\mathbf{g},\mathbf{x}^\top\right]+\left[\mathbf{x},\mathbf{g}^\top\right]}_{B_1} dA^\ast+
		\underbrace{\left[\mathbf{g}^\ast,\mathbf{x}^\top\right]+\left[\mathbf{x},\mathbf{g}^\dagger\right]}_{B_2} dA   .
	\end{dmath*}
	Since both $\mathbf{f}_C$ and $\mathbf{g}$ are commutator-conservative with respect to $\mathbf{x}$, by Proposition \ref{prp:CommutativityPreservationClosed}  the terms $A_1$, $B_1$, $B_2$ are equal to zero. It remains to prove that $A_2 + A_3$ is equal to zero. Since $\mathbf{g}$ is commutator-conservative, then we can write $g_1 = i \pdv*{ \mathbb{L}}{p}$, and $g_2 = -i \pdv*{ \mathbb{L}}{q}$. Using \eqref{eq:fL}, we can expand the terms $A_2$ and $A_3$ as follows,
	\begin{dgroup*}
		\begin{dmath*}
			\left[\mathbf{f}_L,\mathbf{x}^\top\right]+\left[\mathbf{x},\mathbf{f}_L^\top\right] = \begin{bmatrix}
				0 && \frac{\commutator{f_{L,1}}{p}+ \commutator{q}{f_{L,2}}}{2}\\
				\frac{\commutator{f_{L,2}}{q}+ \commutator{p}{f_{L,1}} }{2} && 0		
			\end{bmatrix},
		\end{dmath*}
		\begin{dmath*}
			\left[\mathbf{g}^\ast,\mathbf{g}^\top\right] = \begin{bmatrix}
				0 && g_1^\ast g_2 - g_2^\ast g_1\\
				g_2^\ast g_1 - g_1^\ast g_2 && 0
				
			\end{bmatrix}.
		\end{dmath*}
		
	\end{dgroup*}
	Now we can evaluate the non zero elements of the last two matrices. Substituting for $g_1$ and $g_2$ and using the product rule in Proposition \ref{prp:BasicPropetiesDerivationOnP}\ref{prt:ProductRule}, we obtain after some calculations
	\begin{dmath*}
		\frac{1}{2} \qty(\commutator{f_{L,1}}{p}+ \commutator{q}{f_{L,2}} ) = 
		{\frac{i}{2} \qty( \pdv{ \qty(g_1^\ast \mathbb{L} + \mathbb{L}^\ast g_1)}{q} + \pdv{ \qty(g_2^\ast \mathbb{L} + \mathbb{L}^\ast g_2) }{p}) = g_2^\ast g_1 - g_1^\ast g_2}.
	\end{dmath*}
	Similarly, one can obtain that
$		\frac{1}{2} \qty(\commutator{f_{L,2}}{q}+ \commutator{p}{f_{L,1}}  ) = g_1^\ast g_2 - g_2^\ast g_1 $.
	Hence, we have $A_2 + A_3$ is equal to zero which completes the proof.
\qed\end{pf}

\section{Generalization to Multiple Modes}
In this section, we will generalize the results in the previous section to the case of multiple modes interacting with $n$ environment fields. We assume that the scattering matrix is constant; i.e.,  $\mathbb{S} \in \mathbb{C}^{n \times n}$. In this case, the system dynamics can be described by QSDEs of the form
\begin{dmath}
	d\mathbf{x} = \mathbf{f}(\mathbf{x})dt + \mathbf{g}(\mathbf{x}) d\mathbf{A}^\ast + \mathbf{g}(\mathbf{x})^\ast d\mathbf{A},\label{eq:HeisenBergNonLinearQSDEMulti}
\end{dmath}
From \eqref{eq:QSDE_X}, the functions $\mathbf{f}$ and $\mathbf{g}$ in these QSDEs will be given by,
\begin{dgroup}
	\begin{dmath}
		\mathbf{f}(\mathbf{x}) = -i\left[\mathbf{x},\mathbb{H}\right] + \frac{1}{2}\left[\qty(\mathbb{L}^\dagger \mathbb{S}\mathbf{g}(\mathbf{x})^\top)^\top +  	\qty(\mathbb{L}^\dagger \mathbb{S}\mathbf{g}(\mathbf{x})^\top)^\dagger\right], \label{eq:quantumMarkovGeneratorMulti}
	\end{dmath}
	\begin{dmath}
		\mathbf{g}(\mathbf{x}) =   \left[\mathbf{x},\mathbb{L}^\top\right]\mathbb{S}^\ast. \label{eq:GMulti}
	\end{dmath}
\end{dgroup}
We also define physical realizability as in the previous section as below.
\begin{defn}\label{def:PhysicalRealizabilityMulti}
	QSDEs of the form \eqref{eq:HeisenBergNonLinearQSDEMulti}, where $\mathbf{f} \in  \mathcal{O}_{\mathbf{x}}$ and $\mathbf{g} \in  \mathcal{P}_{\mathbf{x}}$ are said to be physically realizable if there exists a pair $\mathbb{H} \in  \mathcal{O}_{\mathbf{x}},  \mathbb{L} \in  \mathcal{P}_{\mathbf{x}}^{n \times 1}$ and a matrix  $\mathbb{S} \in \mathbb{C}^{n \times n}$ satisfying \eqref{eq:quantumMarkovGeneratorMulti} and \eqref{eq:GMulti}.		
\end{defn}
As in the previous section, we now define the notion of a commutator-conservative mapping for the multiple mode case.
\begin{defn}
	A mapping  $\mathbf{g} \in  \mathcal{P}_{\mathbf{x}}^{2 m \times n}$ is commutator-conservative if there exists $\mathbb{J} \in  \mathcal{P}^{n \times 1}_{\mathbf{x}}$ such that  $\mathbf{g} =  -i\left[\mathbf{x},\mathbb{J}^\top\right]$.
\end{defn}
The following proposition gives a necessary and sufficient condition for a mapping to be commutator-conservative in the case of multiple modes. As in  the corollary of Theorem \ref{thm:ClosedPhysicalRealizability}, it also relates the commutator-conservative mapping and gradient concepts for a multiple channel quantum system. 
Let 
\begin{dmath}
	\bm{\Sigma}_m = \begin{bmatrix}
		0 && \mathbf{I}_m \\ -\mathbf{I}_m && 0
	\end{bmatrix}.
\end{dmath}
\begin{prop}\label{prp:ClosedPhysicalRealizabilityOnPMulti}	
	Consider   $\mathbf{g} \in \mathcal{P}_{\mathbf{x}}^{2 m \times n}$. $\mathbf{g}$ is commutator-conservative
	if and only if $\mathbf{g}$ is a gradient with respect to $\mathbf{y} = \bm{\Sigma}_m \mathbf{x} $.
\end{prop}
\begin{pf}
	First we can decompose $\mathbf{g}$ as
	\begin{dmath*}
		\mathbf{g} = \mqty[\mathbf{g}_1^\top \\ \cdots \\ \mathbf{g}_{2m}^\top], {\mathbf{g}_i \in \mathcal{P}^{n \times 1}_{\mathbf{x}}, \forall i \in \qty{1, \cdots, 2m} }.
	\end{dmath*}
	If $\mathbf{g}$ is a gradient with respect to $\mathbf{y} = \bm{\Sigma}_m \mathbf{x} = \qty[p_{1} \cdots p_{1} \quad \qty(-q_{1}) \cdots \qty(-q_{m})]^\top$, by Theorem \ref{thm:gradient} we can write, for $1\leq i\leq 2m$,
$
		\mathbb{J} = \int_o \mathbf{g}_i dy_i + R(y_i),
$
	where $R(y_i) \in \bar{\mathcal{P}}_{y_i}$. Hence for $1\leq i\leq m , x_i = q_{i}$, we have $-i\commutator{q_{i}}{\mathbb{J}} = \pdv*{\mathbb{J}}{p_{i}} = \pdv*{\mathbb{J}}{y_i} = \mathbf{g}_i$. Furthermore, for $ m+1\leq i\leq 2m ,x_i = p_{i}$, we have $-i\commutator{p_{i}}{\mathbb{J}} = -\pdv*{\mathbb{J}}{q_{i}} = \pdv*{\mathbb{J}}{y_i} = \mathbf{g}_i$. Therefore, $\mathbf{g}$ is commutator-conservative. 
	\\
	Now suppose $\mathbf{g}$ is commutator-conservative.
Then there exists a $\mathbb{J} \in \mathcal{P}_{\mathbf{x}}^{n \times 1}$ such that
	\begin{dmath*}
		{\mathbf{g}
			= -i\left[\mathbf{x},\mathbb{J}^\top\right] = \frac{1}{i} \begin{bmatrix}
				\left[\mathbf{q},\mathbb{J}^\top\right] \\ 
				\left[\mathbf{p},\mathbb{J}^\top\right]
			\end{bmatrix}
			= \frac{1}{i} \begin{bmatrix}
				i\qty(\pdv*{\mathbb{J}}{ \mathbf{p}})^\top \\ -i \qty(\pdv*{\mathbb{J}}{ \mathbf{q}})^\top
		\end{bmatrix}}\\
	= \begin{bmatrix}
		\qty(\pdv*{\mathbb{J}}{ \mathbf{p}})^\top \\ -\qty(\pdv*{\mathbb{J}}{ \mathbf{q}})^\top
	\end{bmatrix}.
\end{dmath*}
Hence $\mathbf{g} = \pdv*{\mathbb{J}}{\mathbf{y}}^\top$ with  $\mathbf{y} = \bm{\Sigma}_m \mathbf{x}$.
\qed\end{pf}

\begin{cor}
	Consider   $\mathbf{g} \in \mathcal{P}_{\mathbf{x}}^{2 m \times n}$. $\mathbf{g}$ is commutator-conservative 
	if and only if
	\begin{dmath}
		\pdv{\mathbf{g}_i}{y_j} = \pdv{\mathbf{g}_j}{y_i}, 
	\end{dmath}
for all $ 1\leq i,j\leq 2m, i\neq j$, and  $\mathbf{y} = \bm{\Sigma}_m \mathbf{x} $.
\end{cor}
\begin{pf}
	This result is immediate from Theorem \ref{thm:Curl}.
\qed\end{pf}

\begin{cor}\label{cor:commutatorConservativeS_ast}
	Let   $\mathbf{g} \in \mathcal{P}_{\mathbf{x}}^{2 m \times n}$ be commutator-conservative. Then for any $\mathbb{S} \in \mathbb{C}^{n\times n}$, such that $\mathbb{S}\mathbb{S}^\dagger = \mathbb{S}^\dagger\mathbb{S} = \mathbf{I}$, there exists $\mathbb{K} \in \mathcal{P}^{n \times 1}_{\mathbf{x}}$ such that  $\mathbf{g} =  -i\left[\mathbf{x},\mathbb{K}^\top\right]\mathbb{S}^\ast$.
\end{cor}
\begin{pf}
	From Proposition \ref{prp:ClosedPhysicalRealizabilityOnPMulti}, there exists $\mathbb{J} \in  \mathcal{P}^{n \times 1}_{\mathbf{x}}$ such that  $\mathbf{g} =  -i\left[\mathbf{x},\mathbb{J}^\top\right]$. Now let $\mathbb{K} = \mathbb{S}\mathbb{J}$. The result immediately follows.
	\qed\end{pf}
\begin{lem}\label{lem:roleOfConstantMulti}
	Consider the QSDEs given in \eqref{eq:HeisenBergNonLinearQSDEMulti}. Let
	\begin{dgroup}
		\begin{dmath}
			\mathbf{Z} = \left. \int_o  \mathbf{g}_1^\top d\mathbf{p} \right|_{ \mathcal{P}_{p}} - \int_o \mathbf{g}_2^\top d\mathbf{q}  , \label{eq:ZMulti}
		\end{dmath}
		\begin{dmath}
					\mathbb{L} = \frac{1}{i} \mathbb{S} \qty(\mathbf{Z} + \mathbf{C}	), \label{eq:LMulti}
		\end{dmath}
		\begin{dmath}
		\mathbf{f}_L = \frac{1}{2}\left[\qty(\mathbb{L}^\dagger \mathbb{S}\mathbf{g}(\mathbf{x})^\top)^\top +  	\qty(\mathbb{L}^\dagger \mathbb{S}\mathbf{g}(\mathbf{x})^\top)^\dagger\right] \label{eq:fLMulti}
		\end{dmath}		
	\end{dgroup}
	where $\mathbf{C} \in \mathbb{C}^{n\times 1}$ is a constant. Suppose $\mathbf{g}$ is commutator-conservative mapping. Let,
	$\mathbb{L}_1 = -i \qty(	Z + \mathbf{C}_1 )$ and $\mathbf{f}_{L,1} = \frac{1}{2}\qty[\qty( \mathbf{g}^\ast \mathbb{L}_1)+\qty( \mathbf{g}^\ast \mathbb{L}_1)^\ast]$, where $C_1$ is a constant. Then for any $\mathbf{C}_1 \in \mathbb{C}^{n\times 1}$, $\mathbf{f} - \mathbf{f}_{L,1}$ is commutator-conservative if $\mathbf{f} - \mathbf{f}_{L}$ is commutator-conservative.
\end{lem}
\begin{pf}
	Immediate generalization from the proof of Lemma \ref{lem:roleOfConstant}\qed.
\end{pf}
\begin{thm}\label{thm:PhysicalRealizabilityOpenMulti}
	Consider a QSDE given as in \eqref{eq:HeisenBergNonLinearQSDEMulti} and an $\mathbb{S} \in \mathbb{C}^{n \times n}$ satisfying, $\mathbb{S}\mathbb{S}^\dagger = \mathbb{S}^\dagger\mathbb{S}=\mathbf{I}$. Using the notation of Lemma \ref{lem:roleOfConstantMulti},  the following conditions are equivalent
	\begin{enumerate}
		\item The QSDEs in \eqref{eq:HeisenBergNonLinearQSDEMulti} are physically realizable,
		\item Both $\mathbf{g}$ and $\mathbf{f} - \mathbf{f}_L$ are commutator-conservative mappings for some $C \in \mathbb{C}$ in \eqref{eq:LMulti}, and the elements of $\mathbf{f}$ are self-adjoint, 
		\item Both $\mathbf{g}$ and $\mathbf{f} - \mathbf{f}_L$ are commutator-conservative mappings for any $C \in \mathbb{C}$ in \eqref{eq:LMulti}, and the elements of $\mathbf{f}$ are self-adjoint. 
	\end{enumerate}
\end{thm}
\begin{pf}
	First, suppose the QSDEs given in \eqref{eq:HeisenBergNonLinearQSDEMulti} are physically realizable.  We observe that the elements of the function $\mathbf{f}$ in \eqref{eq:quantumMarkovGeneratorMulti} are self-adjoint. From \eqref{eq:GMulti}, the commutator-conservativeness of $\mathbf{g}$ follows from $\mathbf{g} = \comm{\mathbf{x}}{\mathbb{L}}\mathbb{S}^\ast =  -i\comm{\mathbf{x}}{\qty(i \mathbb{S}^\dagger \mathbb{L})}$. Therefore, there exists $\mathbf{C} \in \mathbb{C}^{n \times 1}$ such that $\mathbb{L}$ satisfies \eqref{eq:LMulti}. Furthermore from \eqref{eq:quantumMarkovGeneratorMulti}, we observe that using \eqref{eq:fLMulti} we obtain $\mathbf{f} - \mathbf{f}_L = i\comm{\mathbf{x}}{\mathbb{H}}$, which is also commutator-conservative. Therefore, this implies the condition 2.
	\\
	Now suppose condition 2 holds. Lemma \ref{lem:roleOfConstantMulti} shows that condition 2 implies condition 3.
	\\
	Lastly, assume that condition 3 of the theorem is satisfied. 
	Then	from the commutator-conservativeness of $\mathbf{g}$, by Corollary \ref{cor:commutatorConservativeS_ast}, there exists an $\mathbb{L} \in \mathcal{P}_{\mathbf{x}}$ such that $\mathbf{g} = \comm{\mathbf{x}}{\mathbb{L}}\mathbb{S}^\ast$ given by \eqref{eq:LMulti}. Choose any $\mathbf{C} \in \mathbb{C}^{n \times 1}$ for $\mathbb{L}$ in \eqref{eq:LMulti}. Then we have,
	$\mathbb{L}^\dagger \mathbb{S}\mathbf{g}(\mathbf{x})^\top = i \qty(\mathbf{Z} + \mathbf{C})^\dagger \mathbf{g}(\mathbf{x})^\top $.
	Hence by \eqref{eq:fLMulti}, we have that $\mathbf{f}_L$ is independent of $\mathbb{S}$ and $\mathbf{f}_L$ is self-adjoint. Since $\mathbf{f} - \mathbf{f}_L$ is commutator-conservative and self-adjoint, then by Propositions \ref{prp:ClosedPhysicalRealizabilityOnPMulti} and \ref{prp:SelfAdjointNessOfJ}, there exists an $\mathbb{H} \in \mathcal{O}_{\mathbf{x}}$, such that $\mathbf{f} - \mathbf{f}_L = -i\qty[\mathbf{x},\mathbb{H}]$. By the existence of a Hamiltonian $\mathbb{H}$,   $\mathbb{S} \in \mathbb{C}^{n \times n}$ satisfying, $\mathbb{S}\mathbb{S}^\dagger = \mathbb{S}^\dagger\mathbb{S}=\mathbf{I}$, and  a coupling operator $\mathbb{L}$, it follows that both $\mathbf{f}$ and $\mathbf{g}$ satisfy \eqref{eq:quantumMarkovGeneratorMulti} and \eqref{eq:GMulti} respectively, and therefore the QSDEs are physically realizable.\qed\end{pf}	
\begin{rem}
	Theorem \ref{thm:PhysicalRealizabilityOpenMulti} also shows that any physically realizable QSDEs of the form \eqref{eq:HeisenBergNonLinearQSDEMulti} can be realized independently of the choice of $\mathbb{S} \in \mathbb{C}^{n \times n}$ satisfying, $\mathbb{S}\mathbb{S}^\dagger = \mathbb{S}^\dagger\mathbb{S}=\mathbf{I}$. This adds an extra degree of freedom in realizing nonlinear QSDEs.
\end{rem}
\begin{rem}
	Using a similar assertion as in Proposition \ref{prp:PreservationCommutativity}, one can also show that any QSDE of the form \eqref{eq:HeisenBergNonLinearQSDEMulti} which is physically realizable, also preserves the commutation relations.
\end{rem}
\section{Examples}
\subsection{Physical Realizability of Linear QSDEs}
Let $\mathbf{x} = \qty[q \; p]^\top$ and consider single mode linear QSDEs as follows
\begin{equation}
d\mathbf{x} = \mathbf{A}\mathbf{x} dt + \mathbf{B} dA^\ast  + \mathbf{B}^\ast dA. \label{eq:LinearQSDE}
\end{equation}
Without loss of the generality, let $ \mathbf{B} = i \Sigma \mathbf{C} = i\qty[c_2 \;\; -c_1]^\top$, where $c_i \in \mathbb{C}$ and $\mathbf{A} \in \mathbb{R}^{2 \times 2}$. We will consider conditions under which the linear QSDEs \eqref{eq:LinearQSDE} are physically realizable. By Theorem \ref{thm:PhysicalRealizabilityOpen}, we obtain $\mathbb{L} = \mathbf{C}^\top \mathbf{x}$. Then with $\mathbf{g} = \mathbf{B}$ and $\mathbb{L} = \mathbf{C}^\top \mathbf{x}$, using \eqref{eq:fL}, we have, $\mathbf{f}_L = \Im{\qty(c_2^\ast c_1)} \mathbf{x} = -\frac{\gamma}{2}  \mathbf{x}$, with $\gamma = -2\Im{\qty(c_2^\ast c_1)} \in \mathbb{R}$. Then from the commutator-conservativeness of $\mathbf{f}_C = \mathbf{f} - \mathbf{f}_L$, 
we can solve for the Hamiltonian $\mathbb{H}$ as follows,
\begin{dmath*}
\mathbb{H} = \left. \int_o f_{c,1} dp \right|_{ \mathcal{P}_{p}} - \int_o f_{c,2} dq + C_1= 
\frac{\mathbf{A}_{12} p^2}{2} - \frac{\mathbf{A}_{21} q^2}{2} - \qty(\mathbf{A}_{22} + \gamma/2) pq + C_1. 
\end{dmath*}
Select $C_1 = -i/2 \qty(\mathbf{A}_{22} + \gamma/2)$, so that the Hamiltonian is self-adjoint, and given by
\begin{equation*}
\mathbb{H} = 
\frac{\mathbf{A}_{12} p^2}{2} - \frac{\mathbf{A}_{21} q^2}{2} - \qty(\mathbf{A}_{22} + \gamma/2) \qty(\dfrac{qp+pq}{2}) . 
\end{equation*}
However,
\begin{dmath*}
\mathbb{H} = \int_o f_{c,1} dp - \left. \int_o f_{c,2} dq \right|_{ \mathcal{P}_{q}} + C_2= \qty(\mathbf{A}_{11} + \gamma/2) \qty(\frac{qp+pq}{2}) + \frac{\mathbf{A}_{12} p^2}{2} - \frac{\mathbf{A}_{21} q^2}{2}.
\end{dmath*}
Equating the last two integrals, we obtain $\mathbf{A}_{11} = -\mathbf{A}_{22} - \gamma$. The Hamiltonian can then be written as $\mathbb{H} = \frac{1}{4} \mathbf{x}^\top\qty(\bm{\Sigma}^\top \mathbf{A}+\mathbf{A}^\top\bm{\Sigma})\mathbf{x}$.
One can verify that with $\mathbf{A}$ satisfying this condition, the preservation of the commutation relations for the linear QSDEs above 
$
i \qty(\mathbf{A} \bm{\Sigma} + \bm{\Sigma}\mathbf{A}^\top ) + \comm{\mathbf{B}^\ast}{\mathbf{B}^\top} = 0,
$	
also holds. This condition is equivalent to the physical realizability condition for linear quantum systems given in \cite[Theorem 3.4]{James2008}.

\subsection{Synthesis of a Physically Realizable Nonlinear QSDE}
Let $\mathbf{x} = \qty[q \; p]^\top$. Suppose we want to construct a nonlinear coherent quantum controller. Assume that we have determined that the controller should be given by the following nonlinear QSDEs: 
\begin{equation}
d\mathbf{x} = \mathbf{f}\qty(\mathbf{x}) dt + \mathbf{B} dA^\ast  + \mathbf{B} dA^\ast, \label{eq:NonLinearQSDE}
\end{equation}
with $\mathbf{B}$ as in the previous example. In these nonlinear QSDEs, assume that we know $f_1 = q^3$, but $f_2$ is left arbitrary. We would like to determine $f_2$, so that \eqref{eq:NonLinearQSDE} is physically realizable.
\\
As in the previous example, we obtain $\mathbb{L} = \mathbf{C}^\top \mathbf{x}$, as well as $\mathbf{f}_L = \Im{\qty(c_2^\ast c_1)} \mathbf{x} = -\frac{\gamma}{2}  \mathbf{x}$, for some $\gamma \in \mathbb{R}$. 
Then from the  commutator-conservativeness of $\mathbf{f}_C = \mathbf{f} - \mathbf{f}_L$, using \eqref{eq:Curl_pq} we can solve for $f_{c,2} = f_2+ \gamma/2 p$ as bellow
\begin{equation*}
f_{c,2} = - \int_o \pdv{f_{c,1}}{q} dp + R(q),
\end{equation*}
where $R(q)$ to be chosen so that $f_{c,2}$ is self-adjoint. Doing so we obtain,
\begin{dmath*}
f_{c,2} = - \int_o \qty(3q^2+\gamma/2) dp + R(q) = -\qty(3q^2 p + \gamma/2 p) + R(q).
\end{dmath*}
We can choose $R(q) = i3q$ so that $f_{c,2} = -\frac{3}{2}\qty(q^2 p+pq^2) - \gamma/2 p$. Therefore, we obtain 
\begin{align*}
f_2 = -\frac{3}{2}\qty(q^2 p+pq^2) - \gamma p.
\end{align*}
As in the previous example, we can obtain $\mathbb{H} = 1/2 \qty(q^3p+pq^3) - \gamma/2 (qp+pq) + C$, for any $C \in \mathbb{R}$.
Now, we consider the case where $f_1$ is the analytical function $f_1 = \cos(q) = \sum_{i=0}^{\infty} \qty(-1)^i q^{2i}/{\qty(2i)!}$. It is easily verified that $\pdv*{\cos(q)}{q} = -\sin(q)$ and $\pdv*{\sin(q)}{q} = \cos(q)$. Using the same procedure as before, we obtain
\begin{align*}
f_2 = \frac{1}{2}\qty(\sin(q)p+p\sin(q)) - \gamma p.
\end{align*}
 
\section{Conclusions}
In this article, we have derived algebraic necessary and sufficient conditions for a class of nonlinear QSDEs to be physically realizable. 
We have also given two examples which highlight the application of these results.

\bibliographystyle{elsarticle-num}        				  
\bibliography{ReferenceAbbrvBibLatex}           

\appendix
\section{Appendix}
\subsection{Proof Lemma \ref{lem:DerivationQandP}}
\begin{pf}
	We will only prove  $\pdv*{q^m}{q} = m q^{m-1}$, while $\pdv*{p^m}{p} = m p^{m-1}$ can be proven in a similar way. We will establish the proof by induction. Let $n,m \in \mathbb{N}$ be given. Observe that for $m=1$, by direct substitution, we have $\pdv*{q}{q} = \frac{-1}{i}\comm{p}{q} = 1 = m q^{m-1}$. Now we will evaluate the case of $m > 1$, assuming that for $n = m-1$, $\pdv*{q^{n}}{q} = n q^{n-1}$ holds. Then,
	\begin{dmath*}
		{\pdv{q^m}{q} = \pdv{q q^n}{q} = -\frac{1}{i} \commutator{p}{q q^n} = -\frac{1}{i} \commutator{p}{q}q^n - q  \frac{1}{i} \commutator{p}{q^n} }\\
		{= -\frac{1}{i} \commutator{p}{q}q^n + q  \pdv{q^n}{q} = q^n + n q^{n} = m q^{m-1}}.
	\end{dmath*}
	Hence the assertion is proved by induction.
\qed\end{pf}
\subsection {Proof of Proposition \ref{prp:BasicPropetiesDerivationOnP}}
\begin{pf}
	For \ref{prt:CommutativeRule}, let $X \in \mathcal{P}_{q,p}$. Then
	\begin{dmath*}
		{\pdv{X}{q}{p} = \commutator{p}{\commutator{q}{X}} = p q X + X q p  - \qty(q X p + p X q)} = 
		{   q p X + X p q - \qty(q X p + p X q)  = \commutator{q}{\commutator{p}{X}} }
		= \pdv{X}{p}{q}.
	\end{dmath*}
	For \ref{prt:ProductRule} , we can write,
	\begin{dmath*}
		{\pdv{XY}{p}  = \frac{1}{i} \commutator{q}{XY} = \frac{1}{i} \commutator{q}{X}Y + X \frac{1}{i} \commutator{q}{Y} }\\
		{=  \pdv{X}{p} Y + X\pdv{Y}{p}}.
	\end{dmath*}
	The case for $\pdv{XY}{p}$ is similar and hence \ref{prt:ProductRule} holds true. For \ref{prt:SymmetricDerivations}, we can write
	\begin{dmath*}
		{\qty(\pdv{X}{p})^\ast = \qty(\frac{1}{i} \commutator{q}{X})^\ast = -\frac{1}{i} \qty(q X - X q)^\ast}\\
		{ = \frac{1}{i} \commutator{q}{X^\ast} =  \pdv{X^\ast}{p}}.
	\end{dmath*}
	The case for $\qty(\pdv{X}{p})^\ast$ then follows similarly.	
\qed\end{pf}
\subsection{Proof of Lemma \ref{lem:zero_int_unique}}
\begin{pf}
	We will show this result for the $\int_o \cdot dp $ case only. The other case follows similarly. Without loss of the generality, take $X = \phi_j = q^{k_j}p^{l_j}$. Assume that there are $Y,Y' \in \quotient{\mathcal{P}_{q,p}}{\mathcal{P}_{q}}$ such that  $ Y \neq Y'$ and $\int_o X dp = Y$ , $\int_o X dp = Y'$. In addition, we suppose $Y = \frac{1}{l_j + 1} q^{k_j}p^{l_j + 1} \in \quotient{\mathcal{P}_{q,p}}{\mathcal{P}_{q}}$. Then by Lemma \ref{lem:DerivationQandP}, $\pdv*{Y}{p}=X$. Also, since $Y' \in \quotient{\mathcal{P}_{q,p}}{\mathcal{P}_{q}}$, we have,
	\begin{dmath*}
		{X = \pdv{Y'}{p} = \pdv{\qty(Y + \qty(Y' - Y))}{p} = \pdv{Y}{p} + \pdv{\qty(Y' - Y)}{p} }\\
		{= X + \pdv{\qty(Y' - Y)}{p}}
	\end{dmath*}
	Hence $\pdv*{\qty(Y' - Y)}{p} = 0$. But $Y' - Y \in \quotient{\mathcal{P}_{q,p}}{\mathcal{P}_{q}}$, and therefore, $Y' - Y = 0$, which is a contradiction. The general case then follows by linearity.
\qed\end{pf}
\subsection{Proof of Lemma \ref{lem:Py_Expansion}}
\begin{pf}
	First, we observe that we can expand $\mathcal{P}_{\mathbf{x}}$ as follows 
	\begin{dmath}
		\mathcal{P}_{\mathbf{x}} 
		= \qty(\quotient{\mathcal{P}_{\mathbf{x}}}{\bar{\mathcal{P}}_{x_1}} )	\oplus \bar{\mathcal{P}}_{x_1} 
		= \qty(\quotient{\mathcal{P}_{\mathbf{x}}}{\bar{\mathcal{P}}_{x_1}} )	\oplus  \left[\qty(\quotient{\bar{\mathcal{P}}_{x_1}} {\bar{\mathcal{P}}_{\qty(x_1,x_2)}} )\oplus \qty(\quotient{\bar{\mathcal{P}}_{\qty(x_1,x_2)}}{\bar{\mathcal{P}}_{\qty(x_1,x_2,x_3)}}) \oplus \cdots \oplus  \quotient{\mathcal{P}_{x_{2m}}}{\mathbb{C}} \oplus \mathbb{C}\right] 
		= \qty( \quotient {\bar{\mathcal{P}}_{\{\}}}{ \bar{\mathcal{P}}_{x_1}}  )	\oplus  \qty(\quotient{\bar{\mathcal{P}}_{x_1}}{\bar{\mathcal{P}}_{\qty(x_1,x_2)}} ) \oplus \qty(\quotient{\bar{\mathcal{P}}_{\qty(x_1,x_2)}}{\bar{\mathcal{P}}_{\qty(x_1,x_2,x_3)}}) \oplus \cdots \oplus  \qty(\quotient {\bar{\mathcal{P}}_{\qty(x_1,\cdots,x_{2m-1})}}{\bar{\mathcal{P}}_{\mathbf{x}}} ) \oplus \bar{\mathcal{P}}_{\mathbf{x}}. \label{eq:Px_expansion}
	\end{dmath}
	Using this fact, any $f \in \mathcal{P}_{\mathbf{x}} $, can be expanded as follows
	\begin{dmath}
		f = \int_o \pdv{f}{x_1} d x_1 + \left. \int_o \pdv{f}{x_2} d x_2 \right|_{\mathcal{P}_{\qty(x_2,\cdots,x_{2m})}} + \cdots + \left. \int_o \pdv{f}{x_i} d x_i \right|_{\mathcal{P}_{\qty(x_i,\cdots,x_{2m})}} + \cdots + \left. \int_o\pdv{f}{x_{2m}} d x_{2m} \right|_{\mathcal{P}_{x_{2m}}} + C,  \label{eq:f_expandX}
	\end{dmath}
	where $C \in \mathbb{C}$ is a constant. From the fact that $\mathcal{P}_{\mathbf{x}}= \mathcal{P}_{\mathbf{y}}$, we could also expand $	\mathcal{P}_{\mathbf{y}}$ as in \eqref{eq:Px_expansion}. Therefore \eqref{eq:f_expandY} follows.
\qed\end{pf}

\subsection{Proof of Theorem \ref{thm:gradient}}
\begin{pf}
	Necessity. If $f$ can be expanded using elements of $\mathbf{g}$ for any permutation of $\mathbf{x}$, then by Lemma \ref{lem:Py_Expansion}, for any $i$, 	
	$
	f = \int_o g_i  d x_i + R_{x_i},
	$
	where $R_{x_i} \in \bar{\mathcal{P}}_{x_i}$. Since $\pdv{ R_{x_i}}{ x_i} = 0$, it follows immediately that $\pdv*{f}{\mathbf{x}} = \mathbf{g}^\top$.
	\\
	Sufficiency. Let $g_i =  \pdv{f}{x_i}$. Then by \eqref{eq:f_expandY}, it follows immediately that we can write for any $i$ 
	$,		f = \int_o g_i  d x_i + R_{x_i}$
	, where $R_{x_i} \in \bar{\mathcal{P}}_{x_i}$, which completes the proof.
	\qed\end{pf}

\subsection{Proof Theorem \ref{thm:Curl}}
\begin{pf}
	We first prove the sufficiency part. Let $\mathbf{g}$ be a gradient with respect to $\mathbf{x}$. From Theorem \ref{thm:gradient}, and the fact that $\pdv*{x_i}$ and $\pdv*{x_j}$ are commutative, taking two derivatives of the integral equality for both $\pdv{x_i}$ and $\pdv{x_j}$, we obtain
	\begin{dmath*}
	\pdv{}{x_i}{x_j} \qty[\int_o g_i  d x_i + R_{x_i}] = \pdv{}{x_i}{x_j} \qty[\int_o g_j  d x_j + R_{x_j} ],
	\end{dmath*}
	where $R_{x_i} \in \bar{\mathcal{P}}_{x_i}$, and $R_{x_j} \in \bar{\mathcal{P}}_{x_j}$. Consequently we obtain,
	$
	\pdv{g_i}{x_j} = \pdv{g_j}{x_i}.
	$
	\\
	We will establish the necessary part by induction.
	Let $P(n)$ be the statement that if \eqref{eq:CurlGeneralized} is true for all $i \neq j, 1\leq i\leq n, 1\leq j\leq n$, then there exists $f \in \mathcal{P}_{\mathbf{x}}$ such that $\pdv*{f}{\mathbf{x}_i} = \mathbf{g}_i$, for any $1\leq i\leq n$ .
	First we will establish $P(2)$. Suppose \eqref{eq:CurlGeneralized} holds. Then, we can write 
	\begin{subequations}
		\begin{align}
		g_1 = \int_o \pdv{g_2}{x_1}dx_2 + r_1 = 
		\projection{g_1}{\quotient{\mathcal{P}_{\mathbf{x}}}{\bar{\mathcal{P}}_{x_2}}} +\projection{g_1}{\bar{\mathcal{P}}_{x_2}}, & r_1 \in \bar{\mathcal{P}}_{x_2},\\
		g_2 = \int_o \pdv{g_1}{x_2}dx_1 + r_2 = 
		\projection{g_2}{\quotient{\mathcal{P}_{\mathbf{x}}}{\bar{\mathcal{P}}_{x_1}}} +\projection{g_2}{\bar{\mathcal{P}}_{x_1}}, & r_2 \in \bar{\mathcal{P}}_{x_1}.
		\end{align}\label{eq:g1g2_expansion}			
	\end{subequations}
	We observe that by \eqref{eq:Px_expansion}
	\begin{dmath}
	\mathcal{P}_{\mathbf{x}} 
	= \qty(\quotient{\mathcal{P}_{\mathbf{x}}}{\bar{\mathcal{P}}_{x_1}} )	\oplus  \qty(\quotient{\bar{\mathcal{P}}_{x_1}} {\bar{\mathcal{P}}_{\qty(x_1,x_2)}} ) \oplus \bar{\mathcal{P}}_{\qty(x_1,x_2)}. \label{eq:Expansion_PX_version_1}
	\end{dmath} 
	Now according to \eqref{eq:Expansion_PX_version_1}, select $f$ as
	\begin{dmath}
	f = \int_o g_1 dx_1 + \projection{\int_o g_2 dx_2}{\bar{\mathcal{P}}_{x_1}} + r_3, r_3 \in \bar{\mathcal{P}}_{\qty(x_1,x_2)}. \label{eq:f_in_Expansion_version_1}
	\end{dmath}
	Substituting \eqref{eq:g1g2_expansion} to \eqref{eq:f_in_Expansion_version_1}, we obtain,
	\begin{dmath*}
	f = \int_o \int_o \pdv{g_2}{x_1}dx_2 dx_1 + \int_o r_1 dx_1 + \projection{\int_o g_2 dx_2}{\bar{\mathcal{P}}_{x_1}} + r_3.
	\end{dmath*}
	However, according to \eqref{eq:g1g2_expansion} and by the commutativity property of $\int_o dx_1, \int_o dx_2$, we obtain
	\begin{dmath*}
	f = \int_o \int_o \pdv{g_2}{x_1} dx_1 dx_2 + \int_o r_1 dx_1 + \projection{\int_o g_2 dx_2}{\bar{\mathcal{P}}_{x_1}} + r_3.
	\end{dmath*}
	We observe that $\int_o \pdv*{g_2}{x_1}dx_1 = \projection{g_2}{\quotient{\mathcal{P}_{\mathbf{x}}}{\bar{\mathcal{P}}_{x_1}}}$, and , $\int_o r_1 dx_1 + \projection{\int_o g_2 dx_2}{\bar{\mathcal{P}}_{x_1}} = \int_o r_1 dx_1 + \int_o \projection{ g_2 }{\bar{\mathcal{P}}_{x_1}} dx_2$. Therefore, 
	\begin{dmath}
	f = \int_o \projection{g_2}{\quotient{\mathcal{P}_{\mathbf{x}}}{\bar{\mathcal{P}}_{x_1}}} dx_2  + \int_o \projection{ g_2 }{\bar{\mathcal{P}}_{x_1}} dx_2 + \int_o r_1 dx_1 + r_3
	= \int_o g_2 dx_2  + \int_o r_1 dx_1 + r_3 
	= \int_o g_2 dx_2  + \int_o \projection{g_1}{\bar{\mathcal{P}}_{x_2}} dx_1 + r_3 
	= \int_o g_2 dx_2  +  \projection{\int_o g_1 dx_1}{\bar{\mathcal{P}}_{x_2}}  + r_3. \label{eq:f_for_P_2}
	\end{dmath}
	Hence, $P\qty(2)$ is true.
	\\
	Now, we will consider $P\qty(n)$ for $n > 2$ assuming that $P\qty(n-1)$ is true. Since $P(n-1)$ true, then for any permutation $\sigma : (1,\cdots,n-1) \rightarrow (1,\cdots,n-1)$, there exists an $f \in \mathcal{P}_{\mathbf{x}}$ such that
	\begin{dmath}
	f = \int_o \hat{g}_1 d\hat{x}_1 + \projection{\int_o \hat{g}_2 d\hat{x}_2}{\bar{\mathcal{P}}_{\hat{x}_1}} + \cdots + \projection{\int_o \hat{g}_{n-1} d\hat{x}_{n-1}}{\bar{\mathcal{P}}_{\hat{x}_1,\cdots,\hat{x}_{n-2}}} + r_n, r_n \in \bar{\mathcal{P}}_{\qty(x_1,\cdots,x_{n-1})}, \label{eq:f_in_Expansion_version_z}
	\end{dmath}
	where $\hat{x}_i = x_{\sigma(i)}$ and $\hat{g}_i = g_{\sigma(i)}$. Now suppose
	\begin{dmath}
	r_n =  \projection{\int_o g_n dx_n}{\bar{\mathcal{P}}_{\qty(x_1,\cdots,x_{n-1})}} + r_{n+1} \;, r_{n+1} \in \bar{\mathcal{P}}_{\qty(x_1,\cdots,x_n)}. \label{eq:r_n}
	\end{dmath}
	We can write for any $i \leq n-1 $
	\begin{dmath}
	g_i = \int_o \pdv{g_n}{x_i}dx_n + r_i = 
	\projection{g_i}{\quotient{\mathcal{P}_{\mathbf{x}}}{\bar{\mathcal{P}}_{x_n}}} +\projection{g_i}{\bar{\mathcal{P}}_{x_n}}.\label{eq:gngi_expansion}			
	\end{dmath}
	For simplicity, let us assign
	\begin{dmath}
	\int_o \mathbf{\hat{g}}^\top d\mathbf{\hat{x}} = \int_o \hat{g}_1 d\hat{x}_1 + \projection{\int_o \hat{g}_2 d\hat{x}_2}{\bar{\mathcal{P}}_{\hat{x}_1}} + \cdots + \projection{\int_o \hat{g}_{n-1} d\hat{x}_{n-1}}{\bar{\mathcal{P}}_{\hat{x}_1,\cdots,\hat{x}_{n-2}}}.
	\end{dmath}
	Observe that we can also write
	\begin{dmath}
	\mathcal{P}_{\mathbf{x}} = \qty(\quotient{\mathcal{P}_{\mathbf{x}}}{\bar{\mathcal{P}}}_{\hat{x}_1} ) \oplus \qty(\quotient{\bar{\mathcal{P}}_{\hat{x}_1}}{\bar{\mathcal{P}}_{(\hat{x}_1,\hat{x}_2)}} ) \oplus \cdots \qty(\quotient{\bar{\mathcal{P}}_{(\hat{x}_1,\cdots,\hat{x}_{n-2})} }{\bar{\mathcal{P}}_{(\hat{x}_1,\cdots,\hat{x}_{n-1})} } ) 
	\oplus \bar{\mathcal{P}}_{(\hat{x}_1,\cdots,\hat{x}_{n-1})} .
	\end{dmath}
	Substituting \eqref{eq:gngi_expansion} and \eqref{eq:r_n} to \eqref{eq:f_in_Expansion_version_z} we have
	\begin{dmath*}
	f = \int_o \int_o \pdv{g_n}{\hat{x}_1}dx_n d\hat{x}_1 + \projection{\int_o \int_o \pdv{g_n}{\hat{x}_2}dx_n d\hat{x}_2}{\bar{\mathcal{P}}_{\hat{x}_1}} + \cdots + \projection{\int_o \int_o \pdv{g_n}{\hat{x}_{n-1}}dx_n d\hat{x}_{n-1}}{\bar{\mathcal{P}}_{\hat{x}_1,\cdots,\hat{x}_{n-2}}} + \projection{\int_o g_n dx_n}{\bar{\mathcal{P}}_{\qty(x_1,\cdots,x_{n-1})}}\\
	 + \int_o r_1 d\hat{x}_1 + \projection{\int_o r_2 d\hat{x}_2}{\bar{\mathcal{P}}_{\hat{x}_1}} + \cdots + \projection{\int_o r_{n-1} d\hat{x}_{n-1}}{\bar{\mathcal{P}}_{\hat{x}_1,\cdots,\hat{x}_{n-2}}} + \projection{\int_o g_n dx_n}{\bar{\mathcal{P}}_{\qty(x_1,\cdots,x_n)}} 
	+
	r_{n+1} 
	= \int_o \int_o \pdv{g_n}{\hat{x}_1}d\hat{x}_1 dx_n  + \projection{\int_o \int_o \pdv{g_n}{\hat{x}_2}d\hat{x}_2 dx_n }{\bar{\mathcal{P}}_{\hat{x}_1}} + \cdots + \projection{\int_o \int_o \pdv{g_n}{\hat{x}_{n-1}}d\hat{x}_{n-1} dx_n }{\bar{\mathcal{P}}_{\hat{x}_1,\cdots,\hat{x}_{n-2}}} + \projection{\int_o g_n dx_n}{\bar{\mathcal{P}}_{\qty(x_1,\cdots,x_{n-1})}}
	 + \int_o \projection{\hat{g}_1}{\bar{\mathcal{P}}_{x_n}} d\hat{x}_1 + \projection{\int_o \projection{\hat{g}_2}{\bar{\mathcal{P}}_{x_n}} d\hat{x}_2}{\bar{\mathcal{P}}_{\hat{x}_1}} + \cdots + \projection{\int_o \projection{\hat{g}_{n-1}}{\bar{\mathcal{P}}_{x_n}} d\hat{x}_{n-1}}{\bar{\mathcal{P}}_{\hat{x}_1,\cdots,\hat{x}_{n-2}}} + \projection{\int_o g_n dx_n}{\bar{\mathcal{P}}_{\qty(x_1,\cdots,x_n)}} 
	+
	r_{n+1}. 
	\end{dmath*}
	Using the commutativity of zero integration and \eqref{eq:g1g2_expansion}, we obtain
	\begin{dmath*}
	f = \int_o \projection{g_n}{\quotient{\mathcal{P}_{\mathbf{x}}}{\bar{\mathcal{P}}_{\hat{x}_1}} } dx_n  
	+ \projection{\int_o \projection{g_n}{\quotient{\mathcal{P}_{\mathbf{x}}}{\bar{\mathcal{P}}_{\hat{x}_2}} } dx_n }{\bar{\mathcal{P}}_{\hat{x}_1}} 
	+ \cdots 
	+ \projection{\int_o \projection{g_n}{\quotient{\mathcal{P}_{\mathbf{x}}}{\bar{\mathcal{P}}_{\hat{x}_{n-1}}} } dx_n }{\bar{\mathcal{P}}_{\hat{x}_1,\cdots,\hat{x}_{n-2}}} + \projection{\int_o g_n dx_n}{\bar{\mathcal{P}}_{\qty(x_1,\cdots,x_{n-1})}} 	
	 + \projection{ 
		\left[
		\int_o \hat{g}_1 d\hat{x}_1 + \projection{\int_o \hat{g}_2 d\hat{x}_2}{\bar{\mathcal{P}}_{\hat{x}_1}} + \cdots + \projection{\int_o \hat{g}_{n-1} d\hat{x}_{n-1}}{\bar{\mathcal{P}}_{\hat{x}_1,\cdots,\hat{x}_{n-2}}}
		\right]
	}
	{\bar{\mathcal{P}}_{x_n} } 
	+
	r_{n+1}.
	\end{dmath*}
	Therefore, as in \eqref{eq:f_for_P_2}, we can write
	\begin{dmath*}
	f = \int_o \left[\projection{g_n}{\quotient{\mathcal{P}_{\mathbf{x}}}{\bar{\mathcal{P}}_{\hat{x}_1}}}
	+ \projection{g_n}{\quotient{\bar{\mathcal{P}}_{\hat{x}_1}}{\bar{\mathcal{P}}_{(\hat{x}_1,\hat{x}_2)}}}
	+
	\cdots  + \projection{g_n}{\quotient{\bar{\mathcal{P}}_{(\hat{x}_1,\cdots,\hat{x}_{n-2})}}{\bar{\mathcal{P}}_{(\hat{x}_1,\cdots,\hat{x}_{n-1})}}}
	+
	\projection{g_n}{\bar{\mathcal{P}}_{\qty(x_1,\cdots,x_{n-1})}}		
	\right] dx_n  
	+ \projection{\int_o \mathbf{\hat{g}}^\top d\mathbf{\hat{x}}}{\bar{\mathcal{P}}_{x_n}} +
	r_{n+1}\\
	= \int_o g_n dx_n + \projection{\int_o \mathbf{\hat{g}}^\top d\mathbf{\hat{x}}}{\bar{\mathcal{P}}_{x_n}} +
	r_{n+1}.
	\end{dmath*}
		
	Therefore, $P(n)$ is true, which completes the proof by induction.
\qed\end{pf}

\end{document}

%% file: headerSetupAutomatica.tex
\usepackage{amsmath}
\usepackage{amsfonts}
\usepackage{amssymb}
\usepackage{bm}
\usepackage{enumitem}

\usepackage{float}
\usepackage{graphicx}
\usepackage{makeidx}
\usepackage[cal=cm,calscaled=1.2,bb=ams,frak=pxtx,frakscaled=1.2,scr=boondox,scrscaled=1.0,bb=boondox]{mathalfa}

\usepackage{nicefrac}

\usepackage{breqn}

\usepackage{physics}

\usepackage{cite}

\newcommand\quotient[2]{{#1}\setminus{#2}}
\newcommand\projection[2]{\left.#1\right|_{#2}}

\usepackage{etoolbox}

\makeatletter
\preto\maketitle{%
	\begingroup\lccode`~=`,
	\lowercase{\endgroup
		\let\saved@breqn@active@comma~
		\let~}\active@comma 
}
\appto\maketitle{%
	\begingroup\lccode`~=`,
	\lowercase{\endgroup
		\let~}\saved@breqn@active@comma 
}
\makeatother